\documentclass[11pt]{amsart}
\usepackage{latexsym, amssymb, amsmath}
 \DeclareMathOperator{\hd}{Hd}
\DeclareMathOperator{\Hd}{Hd} \DeclareMathOperator{\nd}{Nbhd}
\DeclareMathOperator{\Nd}{Nbhd}

\DeclareMathOperator{\rk}{rank} \DeclareMathOperator{\aut}{Aut}
\DeclareMathOperator{\comm}{Comm}
\DeclareMathOperator{\sl3}{\bold{SL_3}}
\DeclareMathOperator{\pgl3}{\bold{PGL_3}}
\DeclareMathOperator{\sln}{\bold{SL_n}}
\DeclareMathOperator{\pgln}{\bold{PGL_n}}
\DeclareMathOperator{\sp6}{\bold{SP_6}}
\DeclareMathOperator{\pop}{\bold{SO_\Phi}}
\DeclareMathOperator{\psp6}{\bold{PSP_6}}
\DeclareMathOperator{\Aut}{\bold{Aut}}
\DeclareMathOperator{\pl}{\bold{PGL_2}}
\DeclareMathOperator{\G}{\bold{G}}
\DeclareMathOperator{\AG}{\bold{Ad(G)}}
\DeclareMathOperator{\out}{\bold{Out}}
\DeclareMathOperator{\AUT}{Aut} 
\DeclareMathOperator{\h}{\bold{H}} \DeclareMathOperator{\vol}{vol}
\DeclareMathOperator{\st}{\bold{T}}
\DeclareMathOperator{\vrty}{\bold{W}}
\DeclareMathOperator{\del}{\delta }

\begin{document}

\newcommand{\ka}{\kappa}
\newcommand{\ga}{\Gamma}
\newcommand{\tga}{\Gamma}
\newcommand{\tg}{H}
\newcommand{\rt}{\rightarrow}
\newcommand{\gmga}{\tga \backslash \tg}
\newcommand{\neu}{N(\ga)}
\newcommand{\ap}{{\mathcal{A}}}
\newcommand{\lb}{\lambda}
\newcommand{\se}{\subseteq}
\newcommand{\e}{\varepsilon}
\newcommand{\nre}{\sqrt[n]{\varepsilon}}
\newcommand{\tp}{\widetilde{\phi}}
\newcommand{\pp}{\partial \phi}
\newcommand{\app}{{\mathcal{A}}^+}
\newcommand{\ppr}{\partial \phi _R}
\newcommand{\pps}{{\partial \phi _R}_*}
\newcommand{\pu}{U _{\partial}}
\newcommand{\gh}{\Gamma}
\newcommand{\ghp}{\Gamma ^{\varphi}}
\newcommand{\hp}{G^+}
\newcommand{\oppr}{\overline{\partial \phi _R}}
\newcommand{\ta}{K \gamma_\mathfrak{p}}
\newcommand{\de}{\delta}
\newcommand{\tgd}{K \gamma_\mathfrak{p}(\de)}
\newcommand{\p}{\phi}
\newcommand{\pn}{\phi ^{-1}}
\newcommand{\g}{\gamma}
\newcommand{\pgd}{\pi(\gr ,\de)}
\newcommand{\ep}{\hfill $\blacksquare$ \bigskip}
\newcommand{\pru}{\bigskip \noindent {\bf Proof:} \;}
\newcommand{\sa}{\measuredangle _}
\newcommand{\xy}{X_\infty}
\newcommand{\gy}{\gamma_\infty}
\newcommand{\py}{\pi_\infty}
\newcommand{\ey}{e_\infty}
\newcommand{\xr}{X_{\mathfrak{p}}}
\newcommand{\gr}{\gamma_{\mathfrak{p}}}
\newcommand{\pr}{\pi_{\mathfrak{p}}}
\newcommand{\er}{e_{\mathfrak{p}}}
\newcommand{\phr}{\phi_{\mathfrak{p}}}
\newcommand{\phy}{\phi_{\infty}}
\newcommand{\ar}{\G (\mathcal{O}_S)}
\newcommand{\Ar}{\G (\mathcal{O}_S)}

\title{Quasi-isometric rigidity of higher rank $S$-arithmetic lattices}
\author{Kevin Wortman}
\thanks{The author was supported in part by an N.S.F. Postdoctoral Fellowship.}
\date{December 15, 2004}
\thanks{Email: wortman@math.cornell.edu}

\maketitle

\begin{abstract}
We show that $S$-arithmetic lattices in semisimple Lie groups with
no rank one factors are quasi-isometrically rigid.

\end{abstract}

\bigskip \noindent \Large \textbf{1. Introduction} \normalsize
\bigskip

Cocompact lattices in semisimple Lie groups over local fields with
no rank one factors are quasi-isometrically rigid. This was shown
by Kleiner-Leeb \cite{K-L} in general, and Eskin-Farb \cite{E-F 1}
later gave a different proof in the case of real Lie groups.

Eskin then applied the ``quasiflats with holes" theorem for
symmetric spaces of Eskin-Farb \cite{E-F 1} to prove that any
quasi-isometry of a non-cocompact irreducible lattice in a real
semisimple Lie group with no rank one factors is a finite distance
from a commensurator \cite{Es}. As a consequence, any such lattice
is quasi-isometrically rigid. Basic examples of such lattices
include $\sln (\mathbb{Z})$ for $n \geq 3$. Dru\c{t}u has given
another proof of Eskin's theorem \cite{Dr} using asymptotic cones
and the results of \cite{K-L}.

Eskin's theorem has a place in a larger body of work of Schwartz,
Farb-Schwartz, and Eskin. In particular, it has been shown that
any quasi-isometry of an irreducible non-cocompact lattice in a
semisimple real Lie group, which is not locally isomorphic to
$\mathbf{SL_2}(\mathbb{R})$, is a finite distance from a
commensurator (\cite{S}, \cite{F-S}, \cite{S 2}, and \cite{Es});
see \cite{Fa} for a full account.

While the theorem of Kleiner-Leeb applied to cocompact
$S$-arithmetic lattices in semisimple Lie groups with no rank one
factors, the question of quasi-isometric rigidity for
non-cocompact $S$-arithmetic lattices remained unexplored for a
few years. The first account of quasi-isometric rigidity for
non-cocompact $S$-arithmetic lattices (and the only account aside
from this paper) was given by Taback \cite{Ta}. Taback's theorem
states that any quasi-isometry of $\mathbf{SL_2}(\mathbb{Z}[1/p])$
is a finite distance in the sup-norm from a commensurator. Thus,
Taback's theorem provided evidence that quasi-isometries of
$S$-arithmetic lattices could be characterized in the same way as
their arithmetic counterparts.

Following the work of Eskin, we apply the quasiflats with holes
theorem of \cite{W} for products of symmetric spaces and Euclidean
(affine) buildings to show that non-cocompact $S$-arithmetic
lattices in semisimple Lie groups with no rank one factors are
quasi-isometrically rigid. Examples of such lattices include
$\sln(\mathbb{Z}[1/p])$ and $\sln(\mathbb{F}_{q}[t])$ for $n \geq
3$, where $\mathbb{F}_{q}[t]$ is a polynomial ring with
indeterminate $t$ and coefficients in the finite field
$\mathbb{F}_q$. (See Section 5 for more examples.)

As a special case of our results, we show that any finitely
generated group quasi-isometric to $\sln(\mathbb{Z}[1/p])$, is in
fact isomorphic to $\sln(\mathbb{Z}[1/p])$ ``up to finite groups"
as long as $n \geq 3$.

Our proof also shows that cocompact lattices in semisimple Lie
groups with no rank one factors are quasi-isometrically rigid,
thus providing a unified proof of the theorems of Kleiner-Leeb,
Eskin-Farb, and Eskin. In particular, we give a new proof of the
theorem of Kleiner-Leeb---a proof which does not use the theory of
asymptotic cones.

\bigskip \noindent \textbf{Summary of definitions to come.}
In order to state our results, we briefly provide some
definitions. We will expand on these definitions in Section 2.

For any topological group $H$, we let $\aut (H)$ be the
 group of topological group automorphisms of $H$.

 For any valuation $v$ of a global field $K$, let $K_v$ be the
 completion of $K$ with respect to $v$. If $S$ is a set of
 valuations of $K$, then we let $\mathcal{O}_S \leq K$ be the ring of
 $S$-integers.

 We call an algebraic $K$-group $\G$ \emph{placewise not rank one} with respect to $S$ if
 $K_v -\rk (\G)\neq 1$ for all $v \in S$. We denote the adjoint
  representation by $\mathbf{Ad}$, and we let $G$ be the direct product of the groups
 $\AG(K_v)$ over all $v \in S$ for which $\G$ is $K_v$-isotropic.

 Last, we let $\mathcal{QI}(\ar)$ be the quasi-isometry group of $\ar$, and
 we let $\comm (\ar)$ be the commensurator group of
 $\ar$. We warn the reader here that our definition of
 $\comm (\ar)$ is slightly atypical (see Section 2).

\bigskip \noindent \textbf{Quasi-isometries of $S$-arithmetic
 groups.}
 Our main result is

\bigskip \noindent \textbf{Theorem 1.1} \emph{Let $K$ be a global
field and $S$ a finite nonempty set of inequivalent valuations
containing all of the archimedean ones. Suppose $\G$ is a
connected simple $K$-group that is placewise not rank one with
respect to $S$.}

\medskip

\begin{quote}
\noindent \emph{(i)} \emph{If $\G$ is $K$-isotropic and $K$ is a
number field, then there is an isomorphism of topological groups:
$${\mathcal{QI}}\big(\ar\big) \cong \comm(\ar).$$}

\noindent \emph{(ii)} \emph{If $\G$ is $K$-isotropic and $K$ is a
function field, then there exist inclusions of topological
groups:}
$$  \comm(\ar) \hookrightarrow {\mathcal{QI}}\big(\ar\big) \hookrightarrow \aut( G).$$
\emph{Furthermore, the image of ${\mathcal{QI}}(\ar)$ in}
$\aut(G)$ \emph{has measure zero.}

\bigskip

\noindent \emph{(iii)} \emph{If $\G$ is $K$-anisotropic, then
there is an isomorphism of topological groups:}
$${\mathcal{QI}}\big(\ar\big) \cong \aut(G).$$
\end{quote}

\bigskip

As an example of Theorem 1.1(i), we have
$$ \mathcal{QI}(\sl3
(\mathbb{Z}[1/p])) \cong \pgl3 (\mathbb{Q}) \rtimes
\mathbb{Z}/2\mathbb{Z},$$ where the topology on the right side of
the isomorphism is induced by the topology of $\mathbb{Q}$ as the
diagonal subspace of $\mathbb{R} \times \mathbb{Q}_p$. This
example is described in more detail in Section 5, where we also
present five other examples.

We note that the theorem above leaves room for improvement, as the
$K$-isotropic case for function fields is not completely
determined. However, results in this case are still slightly
stronger than they are for the fully resolved $K$-anisotropic
case.

\bigskip \noindent \textbf{Quasi-isometric rigidity.} From Theorem
1.1 we can deduce

\bigskip \noindent \textbf{Corollary 1.2} \emph{Suppose $K$, $S$,
 and $\G$ are as in Theorem 1.1, and suppose that $\G$ is of adjoint type.
  Let $\Lambda$ be a finitely generated group, and assume there
is a quasi-isometry $$\phi : \Lambda \rt \ar.$$}

\begin{quote} \emph{(i) If $\G$ is $K$-isotropic and $K$ is a number field,
then there exists a finite index subgroup $\Lambda _S$ of
$\Lambda$ and a homomorphism $\varphi : \Lambda _S \rt \G
(\mathcal{O}_S)$ with a finite kernel and finite co-image such
that
$$\sup_{\lambda \in \Lambda _S} d\Big(\varphi (\lambda ), \phi (\lambda )\Big)<\infty.$$}

\medskip

\noindent \emph{(ii) If $\G$ is $K$-isotropic and $K$ is a
function field, then there exists a finite group $F$ and an exact
sequence
$$1 \rt F \rt \Lambda \rt \Gamma \rt 1, $$ such that $\Gamma$ is a
non-cocompact lattice in $\aut(G)$.}

\medskip

\noindent \emph{(iii) If $\G$ is $K$-anisotropic, then there
exists a finite group $F$ and an exact sequence $$1 \rt F \rt
\Lambda \rt \Gamma \rt 1, $$ such that $\Gamma$ is a cocompact
lattice in $\aut(G)$.}
\end{quote}

\bigskip

\bigskip \noindent \textbf{Bibliographic note.} We will present
 a proof of Theorem 1.1 that covers all of the
cases above, some of which are well known.

Part (iii) of Theorem 1.1 and Corollary 1.2 was shown by
Kleiner-Leeb \cite{K-L}. Part (iii) was also shown when $K$ is a
number field and $S$ equals the set of archimedean valuations by
Eskin-Farb \cite{E-F 1}. (Note that the theorems in \cite{K-L} and
\cite{E-F 1} are stated in equivalent terms of isometries of
Euclidean buildings and/or symmetric spaces.)

Part (i) of Theorem 1.1 and Corollary 1.2 was shown by Eskin
\cite{Es} with the additional assumption that $S$ equals the set
of archimedean valuations. Dru\c tu has also given a proof of (ii)
assuming $S$ is the set of archimedean valuations \cite{Dr}. The
proof in \cite{Dr} uses results from \cite{K-L}.

Corollary 1.2 follows directly from Theorem 1.1 and, for part (i),
Margulis' superrigidity theorem. The proof of this corollary using
Theorem 1.1 is routine. See, for example, Section 9 of \cite{Es}.

\bigskip \noindent \textbf{Similarities and differences between
 our proof and Eskin's.} The proof of Eskin's theorem involves
  studying the large-scale geometry of
symmetric spaces on which higher rank real semisimple Lie groups
act. Our proof of Theorem 1.1 applies the ``quasiflats with holes"
theorem from \cite{W} (which itself is an extensions of the
quasiflats with holes theorem of Eskin-Farb \cite{E-F 1}) to
extend Eskin's proof by allowing for the presence Euclidean
buildings. (Recall that Euclidean buildings are the natural spaces
acted on by semisimple Lie groups over nonarchimedean local
fields.) We rely on many of Eskin's arguments in using large-scale
geometry to construct a boundary function defined almost
everywhere.

Where our proof differs substantially from Eskin's, is in the way
we complete the boundary function. We are forced to confront this
problem with different methods, since the proof in \cite{Es}
relies on the fact that the Furstenberg boundary of a real
semisimple Lie group is a Euclidean manifold. This is not the case
in general, as the Furstenberg boundary of a semisimple Lie group
over a nonarchimedean local field is a Cantor set. Being unable to
rely as heavily on topological arguments, we turn to algebraic
methods to find a completion. (See Section 4 for an expanded
outline of our proof.)

\bigskip \noindent \textbf{Strong rigidity.} Our main result
can be viewed as a strengthening of strong rigidity.

Recall that the strong rigidity
 theorems---first proved by Mostow and later expanded on greatly by Prasad, Margulis,
and Venkataramana---state that any isomorphism between irreducible
lattices in semisimple Lie groups, which are not locally
isomorphic to $\mathbf{SL_2}(\mathbb{R})$, extends to an
isomorphism of the ambient semisimple group. Thus, the ambient
semisimple group is completely determined by the isomorphism class
of a lattice (\cite{Mo}, \cite{Pr 0}, \cite{Pr}, \cite{Mar}, and
\cite{Ve}).

Our result states that the quasi-isometry class alone of an
$S$-arithmetic lattice meeting the conditions of Theorem 1.1 is
enough to determine the ambient semisimple group.

We note that the proofs of strong rigidity in cases (i) and (ii)
of our main theorem (given by Margulis and Venkataramana
respectively) are rooted in ergodic theory. Our unified proof of
cases (i), (ii), and (iii) is based on the large-scale geometry of
symmetric spaces and Euclidean buildings. As such,  we return to
Mostow's original ideas and present a proof that is of a more
geometric nature than the ergodic theoretical proofs of strong
rigidity.

\bigskip \noindent \textbf{Number fields versus function fields.}
Although our results are not complete in the function field case,
we point out that this is only due to the absence of a
characterization of commensurators which does not exist in the
function field case (see Proposition 7.2).

Throughout the portion of the proof dealing with large-scale
geometry, the function field case allows for significant
simplifications. The simplifications stem from the fact that two
Weyl chambers in a Euclidean building are Hausdorff equivalent if
and only if their intersection contains a Weyl chamber. Of course
this is false for symmetric spaces.

\bigskip \noindent \textbf{Acknowledgements.} I thank my Ph.D.
 thesis advisor, Benson Farb, for giving me the
opportunity to work on this problem and for believing I could
solve it.

Thanks to Alex Eskin for helpful insights and for helping me
discover some mistakes I made along the way.

I am also happy to thank Nimish Shah for showing me how to prove
Proposition 7.2 below, Steven Spallone for explaining numerous
mathematical concepts to me over the past five years, and both Dan
Margalit and Karen Vogtmann for suggestions about the exposition
of this paper.

I would like to acknowledge the University of Chicago for
supporting me as a graduate student while I developed the ideas in
this paper, and Cornell University for the pleasant working
environment given to me while I completed its writing.

Last and most important, I am grateful for Barbara Csima, Benson
Farb, and Dan Margalit; their support, encouragement, and patience
made me into a mathematician.

\bigskip \noindent \Large \textbf{2. Definitions} \normalsize
\bigskip

We will take some time now to be precise with our definitions.

\bigskip
\noindent  \textbf{Quasi-isometries.} For constants $\kappa \geq
1$ and $C \geq 0$, a $(\kappa,C)$ \emph{quasi-isometric embedding}
of a metric space $X$ into a metric space $Y$ is a function $\phi
: X \rightarrow Y$ such that for any $x_{1},x_{2} \in X$:
$$\frac{1}{\kappa} d\big(x_{1},x_{2}\big) -C \leq d\big(\phi(x_{1}),\phi(x_{2})\big)
\leq   \kappa d\big(x_{1},x_{2}\big) +C.$$

We call $\phi$ a $(\kappa,C)$ \emph{quasi-isometry} if $\phi$ is a
$(\kappa,C)$ quasi-isometric embedding and there is a number
$D\geq 0$ such that every point in $Y$ is within distance $D$ of
some point in the image of $X$.

\bigskip
\noindent  \textbf{Quasi-isometry groups.} For a metric space $X$,
we define the relation $\sim$  on the set of functions  $X\rt X$
by $\phi \sim \psi$ if
$$\sup _{x\in X} d\big(\phi(x),\psi(x)\big)<\infty.$$

We form the set of all self-quasi-isometries of $X$, and denote
the quotient space modulo $\sim$ by $\mathcal{QI}(X)$. We call
$\mathcal{QI}(X)$ the \emph{quasi-isometry group} of $X$ as it has
a natural group structure arising from function composition. Note
that if $X$ and $Y$ are quasi-isometric metric spaces, then there
is a natural isomorphism $\mathcal{QI}(X) \cong \mathcal{QI}(Y)$.

In addition to a group structure, we also endow $\mathcal{QI}(X)$
with the quotient of the compact-open topology.

\bigskip \noindent \textbf{Word metrics.} A finitely generated group
$\Gamma$ is naturally equipped with a proper left-invariant
\emph{word metric}. This is the metric obtained by setting the
distance between $\gamma \in \Gamma$ and $1\in \Gamma$ to be the
infimum of the length of all words written in a fixed finite
generating set that represent $\gamma$.

The word metric depends on the choice of finite generating set,
but only up to quasi-isometry. Hence, the group
$\mathcal{QI}(\Gamma )$ is independent of the choice of a finite
generating set for $\Gamma$. The topology on
$\mathcal{QI}(\Gamma)$ is also independent of the choice of a
finite generating set since $\Gamma$ is discrete under all choices
of word metrics.

\bigskip \noindent \textbf{$S$-integers.} Recall that finite
algebraic extensions of either $\mathbb{Q}$ or the field
$\mathbb{F}_p(t)$ of rational functions with indeterminate $t$ and
coefficients in a finite field $\mathbb{F}_p$, are called
\emph{global fields}. If $K$ is a global field then we denote the
set of all inequivalent valuations on $K$ by $V_K$, and we denote
the set of all inequivalent archimedean valuations of $K$ by
$V_K^\infty$.

For any valuation $v\in V_K$, let $K_v$ be the topological
completion of $K$ with respect to $v$. The field $K_v$ is a
locally compact nondiscrete field. Any field satisfying these
topological properties is called a \emph{local field}.

For a finite nonempty set of valuations $S \se V_K$ containing
$V_K^\infty $, we define the ring of \emph{$S$-integers} in $K$ to
be
$$\mathcal{O}_S= \{\, x\in K \mid 1 \geq |x|_v  \mbox{ for all } v \in V_K
- S \,\}.$$

\bigskip \noindent \textbf{Rank.} If a simple algebraic
group $\G$ is defined over a field $L$, we say it is an
$L$\emph{-group}. An $L$-group $\G$ is called \emph{$L$-isotropic}
if $L-\rk (\G)>0$, and called \emph{$L$-anisotropic} otherwise.
(Recall that $L-\rk(\G)$ is the maximum dimension of an algebraic
subgroup of $\G$ which is diagonalizable over $L$.)

For a global field $K$ and a simple $K$-group $\G$, let
$V_K^{\G,a} \se V_K$ be the set of valuations $v$ for which $\G$
is $K_v$-anisotropic. Recall that $v \in V_K^{\G,a}$ is equivalent
to the condition that $\G(K_v)$ is compact.

We define $\G$ to be \emph{placewise not rank one} with respect to
a chosen finite set of valuations $S$, if $K_v-\rk (\G) \neq 1$
for all $v\in S$.

\bigskip \noindent \textbf{$S$-arithmetic groups.} A group is
called \emph{$S$-arithmetic} if it is isomorphic to $\ar$ for some
$K$-group $\G$ and for some finite nonempty set $S \se V_K$
containing $V_K^\infty $.

Throughout the remainder, $\G$ is connected, simple, and placewise
not rank one with respect to $S$. Under these conditions it is
well known that $\ar$ is a finitely generated group, so it admits
a proper word metric.

\bigskip \noindent \textbf{Lattices.} A locally compact group $H$
 supports a Haar measure $\mu$. A discrete subgroup $\Gamma < H$
is called a \emph{lattice} if $H / \Gamma$ has finite volume with
respect to $\mu$. This is necessarily the case if $H / \Gamma $ is
compact. Such lattices are called \emph{cocompact}; they are
called \emph{non-cocompact} otherwise.

We write $\AG$ for the image of $\G$ under the adjoint
representation of $\G$. The adjoint representation has a finite
kernel which equals the center of $\G$.

Define $$G = \prod_{v \in S -V_K^{\G,a}} \AG (K_v).$$ The diagonal
homomorphism of $\ar$ into $G$ has a finite kernel. We write the
image of the diagonal homomorphism as $\ar ^\Delta$. The reduction
theory of Borel, Behr, and Harder established that $\ar ^\Delta$
is a lattice in $G$ and that $\ar ^\Delta$ is cocompact if and
only if $\G$ is $K$-anisotropic.

We point out here that $\ar ^\Delta$ is clearly irreducible as a
lattice in $G$. Recall that a lattice $\Gamma < G$ is
\emph{reducible} if $\Gamma$ contains a finite index subgroup of
the form $\Gamma _1 \times \Gamma _2$ where $$\Gamma _i = \Gamma
\cap \prod _{T_i} \AG (K_v),$$ and $T_1$ and $T_2$ nontrivially
partition $S-V_K^{\G,a}$. Otherwise, $\Gamma$ is
\emph{irreducible}.

Let $\aut (G)$ be the group of all topological group automorphisms
of $G$. Since $G$ has a trivial center, it embeds into $\aut (G)$
via inner automorphisms. Furthermore, $G$ is a closed cocompact
subgroup of $\aut(G)$, so $\ar ^\Delta$ is also a lattice in
$\aut(G)$. Furthermore, $\ar ^\Delta$ is cocompact in $\aut(G)$ if
and only if it is cocompact in $G$.

\bigskip \noindent \textbf{Commensurators.} An automorphism $\psi \in \aut
(G)$ \emph{commensurates} $\ar ^\Delta$ if $\psi(\ar ^\Delta) \cap
\ar ^\Delta$ is a finite index subgroup of both $\psi (\ar ^\Delta
)$ and $\ar ^\Delta$.

Define $\comm(\ar) < \aut(G)$ as the group of automorphisms that
commensurate $\ar ^\Delta$. Notice that $\comm(\ar)$ is different
from the standard definition of the commensurator group of $\ar$
in two ways: we project $\ar$ into $G$, and we do not restrict
ourselves to inner automorphisms.

Let $\aut (K)$ be the group of field automorphisms of $K$. There
is an action of $\aut (K)$ on the set of affine $K$-varieties.
Indeed, if $\vrty$ is an affine $K$-variety, then we let
$\sideset{^{\sigma}}{}{\vrty}$ be the variety obtained by applying
$\sigma$ to the coefficients of the polynomials that define
$\vrty$. We define $\aut(K)_{\G}$ to be the group of automorphisms
$\sigma \in \aut(K)$ such that $\sideset{^\sigma}{}{\G}$ is
$K$-group isomorphic to $\G$.

Since valuations are obtained by embedding $K$ into various local
fields, there is an obvious action of $\aut(K)$ on the set of
valuations $V_K$. We let $\aut (K)_{\G,S}$ be the subgroup of
$\aut(K)_{\G}$ consisting of those $\sigma \in \aut(K) _{\G}$ such
that $\sigma (S-V_K^{\G,a})=S-V_K^{\G,a}$.

The group $\aut (K)$ is finite when $K$ is a global field, so both
$\aut(K)_{\G}$ and $\aut (K) _{\G, S}$ are finite also.

We will see in Section 7 that $\comm(\ar)$ is an extension
$$1 \rt \mathbf{Aut(Ad(G))}(K) \rt \comm(\ar) \rt \aut
(K)_{\G,S} \rt 1, $$ where $\mathbf{Aut(Ad(G))}$ is the $K$-group
of algebraic group automorphisms of $\AG$.

If $\G$ is defined over a subfield of $K$ that is fixed pointwise
by $\aut(K)_{\G,S}$, then the above extension splits. Furthermore,
if $\G$ is $K$-split, then there is a split extension $$1 \rt \AG
(K) \rt \mathbf{Aut(Ad(G))}(K) \rt \mathbf{Out(Ad(G))}(K) \rt 1,$$
where $\mathbf{Out(Ad(G))}$ is the $K$-group of outer
automorphisms of $\AG$ (or alternatively the $K$-group of
automorphisms of the Dynkin diagram of $\AG$).

Combining the two remarks above, we have that if $\G$ is $K$-split
and defined over a subfield of $K$ that is fixed pointwise by
$\aut(K)_{\G,S}$, then
$$\comm(\ar) \cong \Big(\AG (K) \rtimes  \mathbf{Out(Ad(G))}(K)
\Big)\rtimes \aut(K)_{\G,S} .$$

Regardless of whether the extensions defining $\comm(\ar)$ split,
$\comm(\ar)$ contains $\AG(K)$ as a finite index subgroup since
the outer automorphism group of a simple algebraic group is
finite. Therefore, we can define a topology on $\comm(\ar)$ by
assigning the topology on $\AG (K)$ to be the subspace topology
resulting from the diagonal embedding
$$\AG(K) \longrightarrow \prod _{v \in S-V_K^{\G,a}} \AG(K_v).$$

\bigskip \noindent \textbf{Examples.} A reader not familiar with
$S$-arithmetic groups is encouraged at this point to skip ahead to
Section 5 where a series of examples is presented.

\bigskip \noindent \Large \textbf{3. Notes} \normalsize
\bigskip

Now that our definitions are in place, we revisit Theorem 1.1.

\bigskip \noindent \textbf{Remarks on Theorem 1.1(i).} In the $K$-isotropic
case for number fields in Theorem 1.1, the group of $K$-rational
points of $\AG$ is a finite index subgroup of $\mathcal{QI}(\ar)$.
Hence the group operation on $\mathcal{QI}(\ar)$ recovers $K$ and
a finite quotient of $\G$. These are two of the three ingredients
used to create $\ar$. The third ingredient, $S$, cannot in general
be recovered from the quasi-isometry group, but it can be
identified up to an element of the finite group
$\text{Aut}(K)_{\G}$.

Let's briefly make the paragraph above more precise.

Theorem 1.1 states that $\mathcal{QI}(\ar)$ is determined up to a
topological group isomorphism as $\comm(\ar)$. By a theorem of
Borel-Tits (\cite{Bo-T} Cor. 6.7), $\AG (K)^+$ is the minimal
finite index subgroup of $\mathcal{QI}(\ar)$  where $\AG (K)^+$ is
the subgroup of $\AG (K)$ generated by the $K$-points of the
unipotent radicals of the $K$-parabolic subgroups of $\AG$.
Therefore, any topological group isomorphism of
$\mathcal{QI}(\ar)$ induces a topological group isomorphism $$f:
\AG (K)^+ \longrightarrow \AG (K)^+,$$ where we assume the domain
of $f$ has the topology derived from $S$.

Another well known theorem of Borel-Tits (\cite{Bo-T} Theorem (A))
states that $f=\beta \circ \sigma ^0$ where $\sigma \in
\text{Aut}(K)_{\G}$, and
$$\sigma ^0 : \AG (K) \longrightarrow \sideset{^\sigma}{}{\AG}(K)$$ is the
homomorphism defined by applying $\sigma$ to the matrix entries of
$\AG(K)$, and $$\beta :
 \sideset{^\sigma}{}{\AG} \longrightarrow \AG$$ is a $K$-isomorphism of
algebraic groups.

Since $f$ is a homeomorphism,  $\sigma$ is a homeomorphism as
well. Therefore, the topology on the image of $\sigma :K \rt K$ is
given by the set $\sigma S$, since $S$ determines the topology of
the domain of $\sigma$.

Note that if $\sigma \in \text{Aut}(K)_{\G}$ and $\beta
:\sideset{^\sigma}{}{\AG} \rt \AG$ is a $K$-isomorphism of
topological groups, then $\sigma ^0$ restricts to an isomorphism
$\AG (\mathcal{O}_S) \cong \sideset{^\sigma}{}{\AG} ( \mathcal{O}_
{\sigma S})$ and $\beta (\sideset{^\sigma}{}{\AG} ( \mathcal{O}_
{\sigma S}))$ is commensurable with $ \AG (\mathcal{O}_{\sigma
S})$ (see e.g. \cite{Mar} I.3.1.1.iv). Hence, recovering $S$ up to
an element of $\text{Aut}(K)_{\G}$ provides us with enough
information to reconstruct $\ar$ up to finite groups. In light of
this, we could not hope for quasi-isometries to pinpoint $S$ any
more than up to an element of $\text{Aut}(K)_{\G}$.

For clarity, we observe that
$$\sln\big(\mathbb{Z}[i,1/(2+i)]\big) \cong \sln\big(\mathbb{Z}[-i,1/(2-i)]\big)$$
is an example of how the set of valuations can fail to be
identified completely by quasi-isometries since, in this example,
the set cannot even be distinguished by isomorphisms of groups.

\bigskip \noindent \textbf{Remarks on Theorem 1.1(iii).} In the $K$-anisotropic
 case, the simple group $\AG$ is encoded in the
quasi-isometry group, but the global field $K$ is not.

For example, examine the quadratic form
$$\Phi = \sum _{i=1}^5 x_i^2.$$ Let $\pop$ be the special
orthogonal group of $\Phi$, so that $\pop$ is
$\mathbb{Q}$-anisotropic and $\mathbb{Q}(\sqrt{11})$-anisotropic.

There are exactly two elements of
$V_{\mathbb{Q}(\sqrt{11})}^\infty$ --- which we name $v_1^\infty$
and $v_2^\infty$ --- and $\mathbb{Q}(\sqrt{11})_{v_i^\infty}\cong
\mathbb{R}$ for $i=1,2$. If we choose the valuation
$v_{(4+\sqrt{11})} \in V_{\mathbb{Q}(\sqrt{11})}$ defined by the
prime ideal $(4+\sqrt{11}) \se \mathbb{Z}[\sqrt{11}]$, then
$\mathbb{Q}(\sqrt{11})_{v_{(4+\sqrt{11})}}$ is isomorphic to the
field of $5$-adic numbers, $\mathbb{Q}_5$.

Let $S=\{v_1^\infty,v_2^\infty,v_{(4+\sqrt{11})}\}$. By the
theorem of Kleiner-Leeb,
$$\mathcal{QI}\Big(\pop \big(\mathcal{O}_S\big)\Big)
 \cong \pop(\mathbb{Q}_5 ). $$ (That $\pop$ is placewise not rank one with respect to
$S$ follows form the fact that $i\in \mathbb{Q}_5$.)

Next, we take our global field to be $\mathbb{Q}$. We let
$S'=\{v^\infty, v_{(5)}\}$, where $v^\infty$ is the archimedian
valuation on $\mathbb{Q}$ and $v_{(5)}$ is the $5$-adic valuation.
Then Kleiner-Leeb's theorem also gives us
$$\mathcal{QI}\Big(\pop \big(\mathcal{O}_{S'}\big)\Big) \cong
\pop(\mathbb{Q}_5 ). $$ Hence, quasi-isometries could not
distinguish between $\mathbb{Q}$ and $\mathbb{Q}(\sqrt{11})$ in
these two examples.

\bigskip \noindent \textbf{Remarks on Theorem 1.1(ii).} My current level of knowledge
for the general $S$-arithmetic group when
 $K$ is a function field and $\G$ is $K$-isotropic is at an
 intermediate level. In this setting we have stronger results than
 in the $K$-anisotropic case, but less is known than in the number
 field case.

There is some evidence that we should be able to remove the
assumption that $K$ is a number field from part (i) of Theorem
1.1. The number field case itself provides evidence that part (i)
should hold for the function field case, and it has been shown
that Theorem 1.1.(i) holds for $\sln (\mathbb{F}_q[t])$ when $n
\geq 3$ \cite{T-W-W}.

The distinction between number fields and function fields in the
$K$-isotropic case exists because our proof for number fields
takes advantage of Ratner's theorem for unipotent flows \cite{Ra}.
Ratner's theorem is a powerful tool, and it appears to be unknown
in positive characteristic.

Note that, in contrast with lattices in semisimple Lie groups over
$p$-adic number fields, lattices in semisimple Lie groups over
function fields can be non-cocompact. In fact, Harder showed that
if $K$ is a global function field and $\G$ is a simple $K$-group,
then $\G$ can be $K$-anisotropic only if $\G$ is of type
$\mathbb{A}_n$ \cite{Har}. Therefore, resolving the $K$-isotropic
case for function fields has heightened importance.

\bigskip \noindent \textbf{Remarks on Corollary 1.2.} In the remarks
on Theorem 1.1(i) it was pointed out that in the $K$-isotropic
case for a number field $K$, the quasi-isometry group of $\ar$
carries the information needed to reconstruct $\ar$. Hence, an
arbitrary finitely generated group $\Lambda$ that is
quasi-isometric to $\ar$ will also carry the information needed to
reconstruct $\ar$ as $\Lambda $ and $\ar$ will have the same
quasi-isometry groups. This is the content of part (i) of
Corollary 1.2.

Note that (i) states that the only way to deform $\ar$ in the
space of all finitely generated groups without moving it outside
of its initial quasi-isometry class is through algebraic methods.

If we knew that Theorem 1.1(i) held in the function field case,
then Corollary 1.2(i) would apply to the function field case as
well. In particular, case (i) of the above corollary holds when
$\ar$ is replaced by $\sln (\mathbb{F}_q[t])$ for $n\geq 3$.

\bigskip \noindent \textbf{Rigidity for groups with poor finiteness
properties.} Any finitely generated group that was previously
known to be quasi-isometrically rigid contains a finite index
subgroup that is simultaneously complex linear, torsion-free, of
type $F_{\infty}$, and of finite cohomological dimension. Thus,
the final comment in the preceding paragraph displays the first
quasi-isometric rigidity result for a finitely generated group
with poor finiteness properties.

Indeed, it is well known that $\sln (\mathbb{F}_{q}[t])$ is not
virtually torsion free. Hence, $\sln (\mathbb{F}_{q}[t])$ is not
complex linear, and any finite index subgroup has infinite
cohomological dimension. In addition, $\sl3 (\mathbb{F}_{q}[t])$
is known \emph{not} to be finitely presentable (a result of Behr
\cite{Be}), and independent work of Abels and Abramenko shows that
the class of groups of the form $\sln (\mathbb{F}_{q}[t])$ where
$n\geq 3$ contains groups of type $F_{k}$, but not of type
$F_{k+1}$ for all $k \geq 1$ (see \cite{Abl} and \cite{Abr}).
Recall that a group $\pi$ is of type $F_k$ if there exists an
Eilenberg-Mac Lane $K(\pi, 1)$ complex with finite $k$-skeleton,
and $\pi$ is of type $F_{\infty}$ if it is of type $F_{k}$ for all
$k$.

\bigskip \noindent \Large \textbf{4. Outline} \normalsize \bigskip

Our proof of Theorem 1.1 borrows heavily from \cite{Es}.

We proceed by realizing any element of $\mathcal{QI}(\ar)$ as a
quasi-isometric embedding $$\p : N(\Gamma) \longrightarrow X,$$
where $X$ is a product of a symmetric space and a Euclidean
building, and $N(\Gamma) \se X$ is a set (defined in Section 8)
that both contains, and is contained in, a metric neighborhood of
a $\ar $ orbit. The existence of such a quasi-isometric embedding
follows from a theorem of Lubotzky-Mozes-Raghunathan \cite{L-M-R}.

Our goal is to show that $\phi$ is within a finite distance of an
element of $\text{Isom}(X)\cong \aut (G)$.

\bigskip \noindent \textbf{Constructing a boundary function defined a.e.}
In logical order, our proof begins with Section 8. Following
Eskin, we apply basic ergodic theory to show that the generic flat
$F \se X$ has most of its volume contained in $N(\Gamma)$. We
denote this generic collection of flats by $\mathfrak{U}$, and we
note that in general, $\mathfrak{U}$ is a proper subset of the set
of all flats in $X$.

For any flat $F \in \mathfrak{U}$, the quasi-isometric embedding
$\p$ restricts to a quasi-isometric embedding $$\Omega_F'
\longrightarrow X,$$ where $\Omega _F' \se F \cap N (\Gamma)$ is a
suitably large subset of $F$. By precomposing with a closest point
projection, we have maps
$$\p _F : F \longrightarrow X.$$

We analyze the image of these maps using the quasiflats with holes
theorem of \cite{W}, and we use the asymptotic behavior of the
images to construct a function $$\pp : \pu \rt \mathcal{B}(G),$$
where $\mathcal{B}(G)$ is the spherical Tits building for $G$ and
$\pu \se \mathcal{B}(G)$ is a subcomplex that has full measure in
the Furstenberg boundary.

For this task, we mostly defer to the proof in \cite{Es} which
covers the case when $X$ is a symmetric space. Indeed, Eskin's
proof uses the geometry of symmetric spaces mostly to establish a
few foundational lemmas. These lemmas are used to analyze the
behavior at infinity of the quasiflats with holes. We supply the
analogous foundational lemmas for the general space $X$, and then
Eskin's proof applies to the more general setting.

\bigskip \noindent \textbf{Continuity of the boundary
function on neighborhoods of faces.} Section 9 is the final
section of this paper. The first three lemmas of the section are
meant as replacements
  for foundational lemmas
 in \cite{Es}, so that we can apply a proof from \cite{Es} to
 derive a fourth lemma: the restriction of $\pp$ to the simplicial neighborhood of a
 face of a maximal simplex
  in $\pu$ is continuous.

\bigskip \noindent \textbf{Completing the boundary function.} Our
goal is to extend the domain of $\pp$ to all of $\mathcal{B}(G)$.
Then we can use Tits' theorem to show that $\pp$ corresponds to an
element of $\aut (G)$. This step is the content of Section 6.
Despite the fact that this section is the third part of our proof
if it were presented in logical order, it is placed in the early
portion of this paper as it is less technical than material from
Sections 8 and 9, and as it contains material unlike that found in
\cite{Es}.

Eskin's approach to finding an extension of $\pp$, for the case
when $K$ is a number field and $S=V_K^\infty$, was to find a
topological completion of $\pp$. A restriction of $\pp$ to a
co-null subset of the Furstenberg boundary is shown to be
bi-H\"{o}lder. Then $\pp$ can be completed to a domain of
$\mathcal{B} (G)$.

Eskin's argument relied on the fact that the Furstenberg boundary
of a real semisimple Lie group is an analytic manifold and a
topological manifold. In contrast, the Furstenberg boundary of a
semisimple Lie group over a nonarchimedean local field is a Cantor
set. Therefore, our approach is forced to deviate from Eskin's at
this point.

We complete $\pp$ algebraically, using the Borel-Tits
classification of abstract homomorphisms between simple groups. We
restrict $\pp$ to a collection of countably many chambers in
$U_\partial$ (a spherical building for $\G$ over global fields)
and argue that the restriction is induced by an injective
homomorphism of rational points of algebraic groups. The
homomorphism is specified by pairs: isomorphisms of algebraic
groups and inclusions of global fields into local fields. We show
the field inclusions are continuous using the continuity of the
boundary function on simplicial neighborhoods of faces of maximal
simplices. Then we extend the restriction to an automorphism of
$G$ by completing the field inclusions. Finally, we show that the
extension of the restriction is also an extension of $\pp$.

See also \cite{Dr} in the case when $K$ is a number field and
$S=V_K^\infty$ for a more combinatorial approach to this problem.

To conclude Section 6, a result of \cite{Es} is applied to show
that the automorphism of $G$ which corresponds to $\pp$,
stabilizes $\ar$ up to Hausdorff equivalence. We denote the group
of all such automorphisms by $\aut _{\hd}(G;\ar)$. Therefore,
$$\mathcal{QI}(\ar) \cong \aut _{\hd}(G;\ar).$$

\bigskip \noindent \textbf{Automorphisms coarsely preserving
lattices.} If $\G$ is $K$-anisotropic, then $G$ and $\ar$ are
Hausdorff equivalent so $\aut _{\hd}(G;\ar)=\aut(G)$. In Section 7
we show that $\aut _{\hd}(G;\ar)$ is a null subset of $\aut(G)$
otherwise. We also show that $\aut _{\hd}(G;\ar) = \comm (\ar)$
when $G$ is $K$-isotropic and $K$ is a number field.

\bigskip \noindent \Large \textbf{5. Examples} \normalsize \bigskip

This section will be especially useful for geometric group
theorists who are not specialists in $S$-arithmetic lattices.

In this section we present six examples illustrating various
aspects of Theorem 1.1. To focus on previously unknown results,
the examples below will all be for the case that $\G$ is
$K$-isotropic and $S \neq V_K^\infty$.

\bigskip \noindent \textbf{Example (A)} The basic global field
is $\mathbb{Q}$. It supports a countably infinite family of
inequivalent valuations (which we think of as metrics for the
global field): an ``infinite" valuation and an $l$-adic valuation
for every prime number $l$. It is well known that these are the
only valuations supported on $\mathbb{Q}$.

The infinite valuation $v^\infty :\mathbb{Q} \rt \mathbb{R}$ is
obtained by embedding $\mathbb{Q}$ into $\mathbb{C}$ and then
restricting the standard metric on $\mathbb{C}$. Any valuation on
a global field that is obtained through an embedding into
$\mathbb{C}$ is called \emph{archimedean}. By completing
$\mathbb{Q}$ metrically with respect to $v^\infty$ we obtain the
real numbers. In the notation of Section 1, this is written as
$\mathbb{Q}_{v^\infty}=\mathbb{R}$.

The only archimedean valuation on $\mathbb{Q}$ is $v^\infty$, but
there are still the \emph{nonarchimedean} $l$-adic valuations
$v_{(l)}$ for prime numbers $l$. First, we define for any integer
$k$, the natural number $\text{deg}_l(k)$ as the exponent of $l$
occurring in the prime factorization of $k$. Then, we define
$v_{(l)} : \mathbb{Q} \rt \mathbb{R}$ by
$$\Big|\frac{n}{m}\Big|_{v_{(l)}}=\exp\Big(\text{deg}_l(m)-\text{deg}_l(n)\Big).$$
 Hence, the defining feature of the $l$-adic
valuation is that it treats the size of powers of $l$ backwards
from what our intuition is used to from the archimedean valuation.
That is $|l^n|_{v_{(l)}} \to 0$ as $n \to \infty$, and $|1/
l^n|_{v_{(l)}} \to \infty$ as $n \to \infty$.

The $l$-adic valuation on $\mathbb{Q}$ is not complete. If we
complete $\mathbb{Q}$ with respect to $v_{(l)}$, we obtain the
\emph{$l$-adic numbers} $\mathbb{Q}_{v_{(l)}}$ which is written
simply as $\mathbb{Q}_l$. The $l$-adic numbers are locally compact
and totally disconnected.

If we fix a prime number $p$ and let $S=\{v^\infty, v_{(p)} \}$,
then
$$\mathcal{O}_S=\{\, x \in \mathbb{Q} \mid 1 \geq |x|_{v_{(l)}} \mbox{ for all
primes } p\neq l \,\}=\mathbb{Z}[1/p].$$

Because $L-\rk (\sl3 )=2$ for all fields $L$, Theorem 1.1 applies
to $\sl3 (\mathbb{Z}[1/p])$. Since $\mathbb{Q}$ admits no
nontrivial automorphisms, the image of $\sl3$ under the adjoint
representation is $\pgl3$, and transpose-inverse is the only outer
automorphism of $\pgl3$, we have $$\mathcal{QI}(\sl3
(\mathbb{Z}[1/p])) \cong \pgl3 (\mathbb{Q}) \rtimes
\mathbb{Z}/2\mathbb{Z}.$$ Notice that as abstract groups,
$$\mathcal{QI}(\sl3 (\mathbb{Z}[1/p])) \cong\mathcal{QI}(\sl3
(\mathbb{Z}[1/l]))$$ for any primes $p$ and $l$. However this
isomorphism is not topological. Indeed, $\mathcal{QI}(\pgl3
(\mathbb{Z}[1/p]))$ is the quotient of a space of functions so it
has a quotient topology descending from the compact-open topology.
This topology is equivalent to the subspace topology on $\pgl3
(\mathbb{Q})$ inherited from the diagonal embedding
$$\pgl3 (\mathbb{Q}) \rt \pgl3 (\mathbb{R}) \times \pgl3
(\mathbb{Q}_p).$$

With this natural topological structure, the sequence of
quasi-isometry classes given by \[
\begin{pmatrix}
1  & 0 & p^{-n} \\
0 & 1 & 0 \\
0 & 0 & 1
\end{pmatrix}
\]
for $n\in \mathbb{N}$ is discrete in $\mathcal{QI}(\pgl3
(\mathbb{Z}[1/p])) $, but not in $\mathcal{QI}(\pgl3
(\mathbb{Z}[1/l]))$. In particular, $\sl3 (\mathbb{Z}[1/p])$ and
$\sl3 (\mathbb{Z}[1/l])$ are not quasi-isometric if $p\neq l$.

\bigskip \noindent  \textbf{Example (B)} Expanding on the previous
example, we let $P$ be any finite set of prime numbers. Then for
the finite set of valuations $S=\{v^\infty\}\cup \{v_{(p)}\}_{p\in
P}$, the ring $\mathcal{O}_S$ is:  $$\{\, x \in \mathbb{Q} \mid 1
\geq |x|_{v_{(l)}} \mbox{ for all primes } l\notin P
\,\}=\mathbb{Z}[1/m_P],$$ where $m_P=\prod _{p\in P}p$.

Expanding on the previous example in another direction, recall
that for any field $L$, the rank of $\sln$ over $L$ is $n-1$.
Hence, as long as $n \geq 3$ we have $$\mathcal{QI}(\sln
(\mathbb{Z}[1/m_P])) \cong \pgln (\mathbb{Q}) \rtimes
\mathbb{Z}/2\mathbb{Z}.$$

Again we note that $\mathcal{QI}(\sln (\mathbb{Z}[1/m_P]))$ has a
natural topology equivalent to the topology obtained via the
diagonal embedding $$\pgln(\mathbb{Q}) \rt \pgln (\mathbb{R})
\times \prod _{p\in P} \pgln(\mathbb{Q}_p).$$ Hence
$\mathcal{QI}(\sln (\mathbb{Z}[1/m_P]))$ becomes ``more discrete"
as the finite set $P$ grows.

Also notice that the semisimple Lie group
$$\pgln (\mathbb{R})
\times \prod _{p\in P} \pgln(\mathbb{Q}_p)$$ is an index two
subgroup of the topological closure of $\mathcal{QI}(\sln
(\mathbb{Z}[1/m_P]))$. Hence, the quasi-isometry class of $\sln
(\mathbb{Z}[1/m_P])$ identifies the ambient semisimple Lie group
that contains $\sln (\mathbb{Z}[1/m_P])$ as a lattice.

\bigskip \noindent  \textbf{Example (C)} Examine the quadratic form
$$\Phi
=x_1^2+2x_2^2-\sqrt{2}x_3^2+\sum_{i=4}^5(x_i^2-x_{i+2}^2).$$ As
$\Phi$ is defined over $\mathbb{Q}(\sqrt{2})$, the special
orthogonal group $\pop$ is a $\mathbb{Q}(\sqrt{2})$-group.

There are exactly two archimedean valuations supported on
$\mathbb{Q}(\sqrt{2})$. They are obtained from the embeddings
$a+\sqrt{2}b \mapsto a+\sqrt{2}b \in \mathbb{C}$ and $a+\sqrt{2}b
\mapsto a-\sqrt{2}b \in \mathbb{C}$. Call these valuations
$v_1^\infty$ and $v_2^\infty$ respectively, and note that
$\mathbb{Q}(\sqrt{2})_{v_1^\infty}$ and
$\mathbb{Q}(\sqrt{2})_{v_2^\infty}$ are each isomorphic to
$\mathbb{R}$ as topological fields, but each in a different way.

We want to add a nonarchimedean valuation to our example. Since
$3$ does not split as a product of two primes in
$\mathbb{Z}[\sqrt{2}]$, there is a unique extension of the
$3$-adic valuation to $\mathbb{Q}(\sqrt{2})$ (written as
$v_{(3)}$), and $\mathbb{Q}(\sqrt{2})_{v_{(3)}} \cong
\mathbb{Q}_3(\sqrt{2})$.

Let $S=\{ v_1^\infty, v_2^\infty, v_{(3)} \}$. Then
$\mathcal{O}_S=\mathbb{Z}[\sqrt{2},1/3]$. We can apply Theorem 1.1
since the rank of $\pop$ over both $\mathbb{Q}(\sqrt{2})$ and
$\mathbb{Q}(\sqrt{2})_{v_2^\infty}$ is $2$, and the rank of $\pop$
over both $\mathbb{Q}(\sqrt{2})_{v_1^\infty}$ and
$\mathbb{Q}(\sqrt{2})_{v_{(3)}}$ is $3$. (That
$\mathbb{Q}(\sqrt{2})_{v_{(3)}}-\rk (\pop )=3$ follows from the
fact that $\sqrt{-2} \in \mathbb{Q}_3$.)

There is a nontrivial element of $\text{Aut}(
\mathbb{Q}(\sqrt{2}))$. Namely $\sigma $ where $\sigma (a+b
\sqrt{2})=a-b\sqrt{2}$. However, while $\sigma S=S$, there is no
$\mathbb{Q}(\sqrt{2})$-isomorphism between $^\sigma \pop$ and
$\pop$. Indeed, $^\sigma \pop$ and $\pop$ are not even isomorphic
over $\mathbb{R}$ as $^\sigma \Phi$ has signature $(5,2)$ and
$\Phi$ has signature $(4,3)$. Hence, $\AUT (
\mathbb{Q}(\sqrt{2}))_{\G ,S}$ is trivial (as is
 $\textbf{Out} (\pop ) $) so Theorem 1.1 yields
 $$\mathcal{QI}\Big(\pop (\mathbb{Z}[\sqrt{2},1/3])\Big) \cong \pop
 (\mathbb{Q}(\sqrt{2})).$$

\bigskip \noindent  \textbf{Example (D)} The symplectic group $\sp6$
has rank $3$ over any field. For the global field $\mathbb{Q}(i)$,
we take the lone archimedean valuation $v^\infty$ (given by
restricting the standard metric on $\mathbb{C}$) along with the
$(2+i)$-adic and the $(2-i)$-adic valuations to comprise the set
$S$. (Note that $2+i$ and $2-i$ are prime in $\mathbb{Z}[i]$.)

Obviously $\mathbb{Q}(i)_{v^\infty} \cong \mathbb{C}$, and because
$(2+i)(2-i)=5$, both $\mathbb{Q}(i)_{v_{(2+i)}}$ and $
\mathbb{Q}(i)_{v_{(2-i)}}$ are isomorphic to $\mathbb{Q}_5$. Now
$$\mathcal{QI}\Big(\sp6 \big(\mathbb{Z}\big[i,\frac{1}{2+i},
\frac{1}{2-i}\big]\big)\Big) \cong \psp6
 (\mathbb{Q}(i)) \rtimes \mathbb{Z} / 2 \mathbb{Z},$$ where
$\mathbf{PSP_6}$ is the adjoint group of $\sp6$. The nontrivial
element of $\mathbb{Z} / 2 \mathbb{Z}$ represents the automorphism
$\sigma$ of $\mathbb{Q}(i)$ defined by $\sigma (i)
 = -i$. Complex conjugation clearly stabilizes $S$, and
 $\sideset{^\sigma}{}{\psp6}=\mathbf{PSP_6}$ since $\psp6$ is defined over $\mathbb{Q}$.

\bigskip \noindent  \textbf{Example (E)} Let $\mathbb{F}_q$ be the
finite field with $q$ elements, and let $\mathbb{F}_q(t)$ be the
field of rational functions with indeterminate $t$ and
coefficients in $\mathbb{F}_q$. This is the primary example of a
global function field. All other global function fields are finite
algebraic extensions of $\mathbb{F}_q(t)$ in analogy with the role
$\mathbb{Q}$ plays for number fields.

The characteristic of $\mathbb{F}_q(t)$ is nonzero so there are no
embeddings of this field into $\mathbb{C}$ and, hence, no
archimedean valuations.

Examine the valuation of $\mathbb{F}_q(t)$ at infinity,
$v_\infty$, defined on quotients of polynomials by $$
\Big|\frac{p(x)}{q(x)}\Big|_{v_\infty}= \exp
\Big(\text{deg}(p(t))-\text{deg}(q(t))\Big).$$ Note that $v_\infty
$ measures the degree of the pole of a rational function at
$\infty \in \mathbb{P}^1(\overline{\mathbb{F}}_q)$, where
$\overline{\mathbb{F}}_q$ is the algebraic closure of
$\mathbb{F}_q$.

We could define an analogous valuation, $v_p$, for every point $p
\in \mathbb{P}^1(\overline{\mathbb{F}}_q)$. The ring of functions
$f \in \mathbb{F}_q(t)$ for which $|f|_{v_p} \leq 1$ for all $p
\in \mathbb{P}^1(\overline{\mathbb{F}}_q)- \{\infty \}$ are
precisely those rational functions which have no poles in
$\mathbb{P}^1(\overline{\mathbb{F}}_q) - \{\infty \}$.
Equivalently, the ring above is simply the ring of polynomials
with indeterminate $t$. In the notation used in Section 1, we have
$\mathcal{O}_S = \mathbb{F}_q[t]$ for $S=\{ v_\infty \}$.

Completing $\mathbb{F}_q(t)$ with respect to $v_\infty$ produces
the locally compact field of formal Laurent series
$\mathbb{F}_q((t^{-1}))$ with indeterminate $t^{-1}$. Hence, we
have by Theorem 1.1 that
$$\mathcal{QI}\Big(\sln (\mathbb{F}_q[t])\Big) < \Big( \pgln \big(
\mathbb{F}_q((t^{-1})) \big) \rtimes \mathbb{Z}/2\mathbb{Z} \Big)
\rtimes
 \aut \big(\mathbb{F}_q((t^{-1}))\big)$$ for all
$n\geq 3$. We remark that $\aut \big(\mathbb{F}_q ((t^{-1}))\big)$
is profinite and in particular is compact.

It will be shown in \cite{T-W-W} however, that for this particular
example the quasi-isometry group is determined exactly as it is in
the number field case. That is,
$$\mathcal{QI}\big(\sln (\mathbb{F}_q[t])\big) \cong \Big( \pgln \big(
\mathbb{F}_q(t) \big) \rtimes  \mathbb{Z}/2\mathbb{Z} \Big)
\rtimes B,$$ where $B$ is a finite solvable subgroup of $\pl
(\mathbb{F}_q)$. Precisely, $B$ is the group of
$\mathbb{F}_q$-points of $\pl \cong \Aut \mathbf{(}\mathbb{P}^1
\mathbf{)}$ that stabilize our distinguished point $\infty \in
\mathbb{P}^1 (\overline{\mathbb{F}}_q)$.

\bigskip \noindent  \textbf{Example (F)} We give a final example
involving function fields for which I do not at this time know of
a proof that the quasi-isometry group is exactly the subgroup of
$\aut(G)$ consisting commensurators.

Examine the smooth elliptic curve $C$ over $\mathbb{F}_5$ given by
the equation $y^2 =t^3 -t$. The field of $\mathbb{F}_5$-rational
functions on $C$ is $\mathbb{F}_5(t,\sqrt{t^3 -t})$, and it is a
separable extension of $\mathbb{F}_5(t)$.

Note that $(t=2, y=1)$ and $(t=1, y=0)$ define points on $C$ which
we name $p$ and $q$ respectively. We define valuations of
$\mathbb{F}_5(t,\sqrt{t^3 -t})$ with respect to the points $p$ and
$q$ as we did in the previous example, and we let $S=\{v_p ,
v_q\}$. Then $\mathcal{O}_S$ is the ring of regular functions on
$C- \{p,q \}$.

Since $[\mathbb{F}_5(t,\sqrt{t^3 -t}) : \mathbb{F}_5(t)]=2$, and
since the point of $C$ given by $(t=2, y=4)$ and the point $p$
each lie above $2 \in \mathbb{P}^1 (\mathbb{F}_5)$, we know by the
so-called fundamental identity of valuation theory that
$\mathbb{F}_5(t,\sqrt{t^3 -t})_{v_p} \cong \mathbb{F}_5 ((t-2))$.

As the point $q \in C$ is the only point on $C$ with $t=1$ (i.e.
$q$ is a point of ramification) $\sqrt{t^3-t} \notin \mathbb{F}_5
(t)_{w_1} $ where $w_1$ denotes the valuation of $\mathbb{F}_5(t)$
at the point $1 \in \mathbb{P}^1(\overline{\mathbb{F}_5})$. Hence,
$\mathbb{F}_5(t,\sqrt{t^3 -t})_{v_q} \cong
\mathbb{F}_5((t-1))(\sqrt{t^3 - t})$.

Now we are set to apply Theorem 1.1 which states in this case that
$$\mathcal{QI}\Big(\sp6 (\mathcal{O}_S)\Big) $$
is contained as a measure zero subgroup of the direct product of
$$ \psp6 \Big(\mathbb{F}_5 ((t-2))\Big) \rtimes
\aut\Big(\mathbb{F}_5 ((t-2))\Big)$$ with
$$ \psp6 \Big(\mathbb{F}_5
((t-1)) (\sqrt{t^3 -t})\Big) \rtimes \aut\Big(\mathbb{F}_5
((t-1))(\sqrt{t^3-t})\Big).$$

This is a stronger result than the one that is known to hold in
the $K$-anisotropic case, but it is an incomplete result. There is
evidence to suggest that there should be an isomorphism
$$\mathcal{QI}\Big(\sp6 (\mathcal{O}_S)\Big) \cong \psp6
\big(\mathbb{F}_5 (t, \sqrt{t^3-t})\big). $$ Note that it can be
shown that $\aut (\mathbb{F}_5 (t, \sqrt{t^3-t}))_{\G,S}$ is
trivial since there are no nontrivial automorphisms of $C$ which
fix the point $p$ and the point $q$.

Corollary 1.2.(i) would hold for $\sp6 (\mathcal{O}_S)$ if the
above isomorphism existed.

\bigskip \noindent \Large \textbf{6. Completing the boundary function}
\normalsize \bigskip

Let $\G (\mathcal{O}_S)$ be as in Theorem 1.1. Since $\G
(\mathcal{O}_S)$ and $\AG (\mathcal{O}_S)$ are commensurable up to
finite kernels (see e.g. \cite{Mar} I.3.1.1.iv),
$$\mathcal{QI}(\G(\mathcal{O}_S))\cong
\mathcal{QI}(\AG(\mathcal{O}_S).$$ Thus we may, and will, assume
throughout the remainder that $\G$ is of adjoint type.

Let $$G= \prod _{v\in S-V_K^{\G ,a}} \G (K_v).$$

Let $X$ be the natural product of nonpositively curved symmetric
spaces and Euclidean buildings on which $G$ acts by isometries and
such that $\text{Isom}(X)/G$ is compact. In this case
$\text{Isom}(X)\cong \aut (G)$.

Throughout we let $n$ equal the rank of $X$. (Recall the rank of
$X$ is the maximal dimension of a flat in $X$.)

\bigskip \noindent \textbf{Two boundaries.} For any point $e \in X$, there is
a natural topology on the space of directions from $e$ which forms
a simplicial complex $\mathcal{B} (G)$, called the \emph{spherical
Tits building for} $G$. The spherical building is
$(n-1)$-dimensional, and it is the same as the spherical building
for $G$ that is produced using the standard BN pair construction.
Hence, group automorphisms of $G$ induce simplicial automorphisms
of $\mathcal{B} (G)$.

A subset $L\se X$ is called a \emph{wall} if it is a codimension
$1$ affine subspace of a flat that is contained in at least two
distinct flats. A \emph{Weyl chamber} in $X$ is the closure of a
connected component of a flat $F \se X$ less all the walls
containing a fixed point $x \in F$. Most of the time we will not
care about the point $x$ which was used to create a Weyl chamber.
In those cases when the distinction is important, we say any such
Weyl chamber is \emph{based} at $x$. (This is different
terminology than was used in \cite{W}. See the word of caution
following the discussion of the Furstenberg metric.)

The \emph{Furstenberg boundary of} $X$ is the compact space of
maximal simplices in $\mathcal{B} (G)$. We denote it by
$\widehat{X}$. It can be defined equivalently as the space of Weyl
chambers in $X$ modulo the relation that two Weyl chambers are
equivalent if they are a finite Hausdorff distance from each
other.

If $X=\xy \times \xr$, where $\xy$ is a symmetric space and $\xr$
and a Euclidean building, then $\widehat{X}=\widehat{X}_\infty
\times \widehat{X}_\mathfrak{p}$.

\bigskip \noindent \textbf{Furstenberg metric.} There are metrics on $\widehat{X}_\infty$ and
$\widehat{X}_\mathfrak{p}$ that are invariant under a fixed
isotropy subgroup of $\text{Isom}(\xy)$ and $\text{Isom}(\xr)$
respectively. The metric on $\widehat{X}_\infty$ is well-known.

To define the metric on $\widehat{X}_\mathfrak{p}$, we begin by
choosing a point $x\in \xr$ and a representative Weyl chamber
$\mathfrak{S} \se \xr$ for every equivalence class in
$\widehat{X}_\mathfrak{p}$ such that $\mathfrak{S}$ is based at
$x$. Thus, we regard $\widehat{X}_\mathfrak{p}$ as the space of
all Weyl chambers based at $x$.

For any Weyl chamber based at $x$, say $\mathfrak{S}$, let $\g
_\mathfrak{S} :[0, \infty ) \rt \mathfrak{S}$ be the geodesic ray
such that $\g _\mathfrak{S}(0) =x$ and such that $\g _\mathfrak{S}
(\infty )$ is the center of mass of the boundary at infinity of
$\mathfrak{S}$ with its usual spherical metric.

We endow $\widehat{X}_\mathfrak{p}$ with the metric
$\widehat{d}_\mathfrak{p}$ where
\begin{equation*}
\widehat{d}_\mathfrak{p} (\mathfrak{Y},\mathfrak{Z}) =
\begin{cases}
\, \quad \pi, &\text{if $\g _\mathfrak{Y} \cap \g _\mathfrak{Z} = \{ x \}$;}\\
\exp \big(-|\g _\mathfrak{Y} \cap \g _\mathfrak{Z} |\big),
&\text{otherwise.}
\end{cases}
\end{equation*}
In the above, $|\g _\mathfrak{Y} \cap \g _\mathfrak{Z} |$ is the
length of the geodesic segment $\g _\mathfrak{Y} \cap \g
_\mathfrak{Z} $.

Note that $\widehat{d}_\mathfrak{p}$ is invariant under the action
of the stabilizer of $x$ and is a complete ultrametric on
$\widehat{X}_\mathfrak{p}$. That $\widehat{d} _\mathfrak{p}$ is an
ultrametric means that it is a metric and
$$\widehat{d}_\mathfrak{p} (\mathfrak{Y},\mathfrak{Z}) \leq \max
\{ \widehat{d}_\mathfrak{p} (\mathfrak{Y},\mathfrak{X}) ,
\widehat{d}_\mathfrak{p}(\mathfrak{X},\mathfrak{Z}) \} \text{
\quad for any }\mathfrak{Y,Z,X} \in \widehat{X}_\mathfrak{p}.$$

We endow $\widehat{X}$ with the metric $\widehat{d}
=\max\{\widehat{d}_\infty , \widehat{d}_\mathfrak{p}\}$.

\bigskip \noindent \textbf{Caution.} In \cite{W}, Weyl chambers in
buildings are called sectors, and the metric
$\widehat{d}_\mathfrak{p}$ is given a different form. In \cite{W},
we made arguments by projecting onto the factors of $X$, and most
of the paper analyzed the geometry of Euclidean buildings. Thus,
our proof was geared towards terminology and tools more common for
buildings. In this paper, we favor terminology and metrics for
buildings which are more compatible with their better established
symmetric space counterparts.

\bigskip \noindent \textbf{A boundary function defined a.e.} In
Section 8, we will define a group $\Gamma$ that acts on $X$ and is
isomorphic to $\ar$ up to finite groups ($\Gamma$ is a lattice in
the simply connected cover of $G$). We will also define a
$\Gamma$-invariant set $N(\Gamma)\se X$ such that $\Gamma
\backslash N(\Gamma)$ is compact. A theorem of
Lubotzky-Mozes-Raghunathan \cite{L-M-R} states that $\Gamma$ is
quasi-isometric to any metric neighborhood of an orbit of $\Gamma$
in $X$. Hence, if we are given a quasi-isometry of $\ar$, we may
replace it with an equivalent quasi-isometric embedding
$$\phi :N(\Gamma) \longrightarrow N(\Gamma
) \se X.$$

Every direction in $X$ (i.e. every geodesic ray) is contained in a
flat. In Section 8 we will show that enough flats in $X$ have
enough of their volume contained in $N(\Gamma)$ to enable us to
construct a boundary function $$\pp : U_\partial \rt \mathcal{B}
(G),$$ where $U_\partial $ is a subcomplex of $\mathcal{B} (G)$
that has full measure in $\widehat{X}$. The function $\pp$ is a
simplicial isomorphism of $U_\partial$ onto its image.

We state below a lemma on a topological property of $\pp$ that is
proved in Section 9. First, we define $\mathcal{N}(f)$ as the
simplicial neighborhood in $\mathcal{B}(G)$ of a fixed
$(n-2)$-dimensional simplex $f \subset \mathcal{B}(G)$. That is,
$\mathcal{N}(f)$ is the set of all chambers in $\mathcal{B} (G)$
containing $f$. We define $\mathcal{N}_U(f)$ to be the simplicial
neighborhood of $f$ in $\pu$, or $\mathcal{N}(f)\cap \pu$.

\bigskip \noindent \textbf{Lemma 9.4} \emph{If $f \subset \pu$ is a simplex of
dimension $n-2$, then $\pp |_{{\mathcal{N}}_U(f)}$ is continuous
in the Furstenberg metric.}

\bigskip

Our goal is to show that $\pp$ is the restriction of an
automorphism of $\mathcal{B} (G)$ which is continuous on
$\widehat{X}$. Then by Tits' Theorem, $\pp$ is induced by an
element of $\aut(G) \cong \text{Isom}(X)$. Knowing this would
enable us to apply an argument of Eskin's to show further that
$\pp$ corresponds to an isometry of $X$ which is a finite distance
in the sup norm from $\phi$.

\bigskip \noindent \textbf{Embeddings of spherical buildings.} An
\emph{embedding of spherical buildings} $\mathcal{B} _1$ into
$\mathcal{B} _2$ is a function $f: \mathcal{B} _1 \rt \mathcal{B}
_2 $ that restricts to a simplicial isomorphism between
$\mathcal{B} _1$ and $f (\mathcal{B} _1)$.

We wish to describe a particularly nice class of embeddings that
play a key role in our proof. These are embeddings which arise
from extremely well behaved homomorphisms of rational points of
simple groups. We begin by describing the latter.

Let $k$ be an arbitrary field and $\h$ a simple $k$-group. If $k'$
is an extension of $k$, then there are injective group
homomorphisms of $\mathbf{H}(k)$ into $\mathbf{H}(k')$ of the form
$\beta \circ \psi ^0$, where $\psi : k \rt k'$ is an injective
homomorphism of fields and $\beta: \sideset{^{\psi}}{}{\h}\rt
\mathbf{H}$ is a $k'$-isomorphism of algebraic groups. Any such
homomorphism will be called \emph{standard}.

Now let $\mathcal{B} ({\bold{H}}(k))$ and $\mathcal{B}
({\bold{H}}(k'))$ be the spherical buildings for $\mathbf{H}(k)$
and $\mathbf{H}(k')$ respectively. A standard homomorphism induces
an embedding $f: \mathcal{B} ({\bold{H}}(k)) \rt \mathcal{B}
({\bold{H}}(k'))$. We call any such embedding \emph{standard} as
well.

Implicit in theorems of Tits and Borel-Tits, is

\bigskip \noindent \textbf{Proposition 6.1} \emph{Let $\bold{H}$ be a
simple connected $k$-group of adjoint type and assume $k$ is
infinite. If $k'$ is an extension of $k$ with
  $k-\rk ( {\bold{H}})=k'-\rk ( {\bold{H}}) \geq 2$,
then any embedding  $\rho: \mathcal{B} ({\bold{H}}(k)) \rt
\mathcal{B} ({\bold{H}}(k'))$ is standard.}

\pru Let $\mathbf{H}(k)^+$ be the subgroup of $\mathbf{H}(k)$
generated by the $k$-points of the unipotent radicals of
$k$-parabolic subgroups of $\mathbf{H}$. In Chapter 5 of
\cite{Ti}, Tits shows how to construct an injective group
homomorphism $\rho _* :{\bold{H}}(k)^+ \rt {\bold{H}}(k')$ which
is induced by $\rho$. We have used the equal rank condition here.

We would like to be able to apply the well known theorem of
Borel-Tits that classifies certain abstract homomorphisms between
rational points of simple groups as being standard (\cite{Bo-T}
Theorem (A)).

By construction, $\rho_*$ has a nontrivial image. Hence, our
assumptions on $\h$ and $k$ satisfy all of the hypotheses on $\rho
_*$ needed to apply the theorem of Borel-Tits except, possibly,
for the condition that the image of $\rho_*$ is Zariski dense in
$\mathbf{H}$. If we let $\bold{M}$ be the the Zariski closure of
the image of $\rho _*$, then our goal is to show that
$\mathbf{M}=\mathbf{H}$.

By Corollary 6.7 of \cite{Bo-T}, we know that ${\bold{H}}(k)^+$
has no proper finite index subgroup. Hence, $\mathbf{M}$ must be
connected. Also note that $\mathbf{M}$ modulo its radical,
$\mathbf{R(M)}$, has positive dimension since ${\bold{H}}(k)^+$ is
not solvable. In particular there exists a connected simple factor
$\mathbf{L}$ of positive dimension of $\mathbf{M/R(M)}$.

We postcompose $\rho _*$ with the natural sequence of
homomorphisms,
$$\mathbf{M}\rt \mathbf{M/R(M)} \rt \mathbf{L} \rt
\mathbf{Ad(L)},$$ to obtain a homomorphism ${\bold{H}}(k)^+ \rt
\textbf{Ad(L)}(k')$ with a nontrivial, Zariski dense image. Now we
can apply Theorem (A) of \cite{Bo-T} to conclude that there exists
a field homomorphism $\psi :k \rt k'$ and an isogeny
$\sideset{^{\psi}}{}{\h} \rt \mathbf{Ad(L)}$. Therefore,
$$\text{dim}(\mathbf{H})=\text{dim}(\sideset{^{\psi}}{}{\h})=\text{dim}
(\mathbf{Ad(L)})\leq \text{dim} (\mathbf{M/R(M)}) \leq
\text{dim}(\mathbf{M}).$$ Because $\mathbf{H}$ is connected and
$\mathbf{M}\leq \mathbf{H}$, we conclude that
$\mathbf{M}=\mathbf{H}$ as desired. We are then able to apply
Theorem (A) of \cite{Bo-T} to our original homomorphism $\rho _*$
and arrive at our desired conclusion.

 \ep

\bigskip \noindent \textbf{A global sub-building.} We would like to
 be able to apply Proposition 6.1 to an algebraically defined
  sub-building of $\mathcal{B}(G)$. We will need to begin by
  finding an
  extension of $K$, for each $v \in S-V_K^{\G,a}$, that is contained in $K_v$ and
   that satisfies the hypothesis
   of Proposition 6.1. This is the purpose of the following

\bigskip \noindent \textbf{Lemma 6.2} For each $v \in S$, there is
a finite algebraic extension $L^v$ of $K$ such that $L^v$ is
contained in $K_v$ and $L^v-\rk (\G)=K_v-\rk(\G)$.

\pru Given a maximal $K_v$-torus $\st<\G$, there is a group
element $g \in \G(K_v)$ such that $\sideset{^{g}}{}{\st}$ is
defined over $K$, where $\sideset{^{g}}{}{\st}$ denotes the
conjugate of $\st$ by $g$. See Section 7.1 Corollary 3 in
\cite{Pl-Ra} for a proof of this fact. It is assumed that $K$ is a
number field throughout most of \cite{Pl-Ra}, but the proof of
this fact does not make an essential use of the number field
assumption, aside from the proof of the $K$-rationality of the
maximal toric variety of $\G$. For a proof of this last fact over
arbitrary fields $K$, see \cite{Bo-Sp}.

Assume that $\st$ and $g$ are as above and that
$K_v-\rk(\st)=K_v-\rk(\G)$. It is well known that there is a
finite separable extension $F^v$ of $K$ over which
$\sideset{^{g}}{}{\st}$ splits (see e.g. \cite{Bo 2} 8.11). Hence,
if $\mathbf{X}(\sideset{^{g}}{}{\st})_L$ is the group of
characters of $\sideset{^{g}}{}{\st}$ defined over an extension
$L$ of $K$, we have $$ \mathbf{X}(\sideset{^{g}}{}{\st})_{K_v}
=\mathbf{X}(\sideset{^{g}}{}{\st})_{F_v} \cap
\mathbf{X}(\sideset{^{g}}{}{\st})_{K_v}=\mathbf{X}(\sideset{^{g}}{}{\st})_{F_v\cap
K_v}.$$ (Recall that a torus splits over a field $L$ if and only
if all of its characters are defined over $L$.)

 Therefore, we
let $L^v=F^v \cap K_v$ so that $$K_v
-\rk(\st)=K_v-\rk(\sideset{^{g}}{}{\st})=L^v-\rk(\sideset{^{g}}{}{\st}).$$
Hence, $$K_v-\rk(\G) \leq L^v-\rk(\G).$$ Since $L^v < K_v$, the
inequality is an equality.

\ep

We define the group
$$G_R=\prod _{v\in S-V_K^{\G,a}} {\bold{G}}(L^v).$$
Let $\mathcal{B} (G_R)$ be the spherical building for $G_R$. By
our choice of $L^v$, the building $\mathcal{B} (G_R)$ has
countably many chambers, the dimensions of $\mathcal{B} (G_R)$ and
$\mathcal{B} (G)$ are equal, and $\mathcal{B} (G)$ naturally
contains $\mathcal{B} (G_R)$ as a subcomplex.

By conjugating $\ar$, we can assume that $\mathcal{B} (G_R) \se
\pu$.
 Indeed, since $\mathcal{B} (G_R)$ has countably many chambers, we can
 appeal to Lemma 8.9 below.

\bigskip \noindent \textbf{Extending the global embedding.} Define
$\ppr$ as the restriction of $\pp$ to $\mathcal{B} (G_R)$. The
induced group homomorphism $$\pp_{R*} :\prod _{v\in S-V_K^{\G,a}}
{\bold{G}}(L^v)^+ \longrightarrow G$$ has a nontrivial image in
each factor of $G$ by construction. Also, Tits proved that each
$\G(L^v)^+$ is an abstract simple group (\cite{Ti 1} Main
Theorem). It follows that $\pp_{R*}$, and hence $\pp_R$, preserves
factors up to permutation.

Therefore we can apply Proposition 6.1 to conclude that $\ppr$ is
induced by a family of standard homomorphisms. Precisely, there is
a permutation $\tau$ of $S-V_K^{\G,a}$, and for each $v \in
S-V_K^{\G,a}$ there exists an injective field homomorphism
$$\psi _v : L^v \rt K_{\tau(v)}$$
and a $K_{\tau (v)}$-isomorphism of algebraic groups
$$\beta _v: \sideset{^{\psi _v}}{}{\G} \rt \G $$
such that $\pp_{R*}$ is the product of the homomorphisms
$$\beta _v \circ \psi _v ^0 :\G(L^v)^+ \rt \G(K_{\tau(v)}).$$

Now extending $\ppr$ amounts to extending each $\psi _v$. This is
the technique of the proposition below. Before we continue though,
we require an extra piece of notation.

Let $f\subset \mathcal{B}(G_R)$ be an $(n-2)$-dimensional simplex.
We denote the simplicial neighborhood of $f$ in
$\mathcal{B}(G_R)$, or $\mathcal{N}(f) \cap \mathcal{B}(G_R)$, by
$\mathcal{N}_R(f)$.

We continue with

\bigskip \noindent \textbf{Proposition 6.3} \emph{The map $\ppr :
 \mathcal{B} (G_R) \rt \mathcal{B} (G)$
uniquely extends to an embedding $\overline {\ppr} : \mathcal{B}
(G) \rt \mathcal{B} (G) $ which is uniformly continuous on the
Furstenberg boundary.}

\pru Choose an apartment $\Sigma \se \mathcal{B}(G_R) \se
\mathcal{B}(G)$ and a chamber $c \se \Sigma$. For any
$(n-2)$-dimensional simplex $f \se c$, there exists a root space
$R_f \se \Sigma$ (as defined in \cite{Ti} 1.12) such that $f \se
\partial R_f$.

By Proposition 3.27 in \cite{Ti}, any chamber in
$\mathcal{N}_R(f)$ is contained in an apartment of
$\mathcal{B}(G_R)$ which contains $R_f$. Therefore, by Proposition
5.6(i) of \cite{Ti}, there exists a valuation $w(f) \in S
-V_K^{\G,a}$ and an $L^{w(f)}$-defined root subgroup
$\mathbf{U_{w(f)}}<\G$, such that $\mathbf{U_{w(f)}}(L^{w(f)})$
acts faithfully and transitively on $\mathcal{N}_R(f)-\{\,c\,\}$.

The valuation $w(f)$ depends on a choice of $f$. However, for any
valuation $v \in V_K^{\G,a}$, we can choose a face $f_v \se c$
such that $w(f_v)=v$. We assume we have chosen such a face $f_v$
for all $v \in S - V_K^{\G,a}$.

If $b_v \in \mathcal{N}_R(f_v)-\{\,c_v\,\}$, then for any $u \in
\mathbf{U_v} (L^v)$ we have $ub_v \in \mathcal{B}(G_R)$.
Therefore,
$$\pp(ub_v)=\ppr(ub_v)=\beta \circ \psi_v^0(u)\ppr(b_v) $$
Since $\mathcal{N}_R(f_v) \se \mathcal{N}_U(f)$, it follows from
Lemma 9.4 that $\beta \circ \psi_v^0$, and hence $\psi_v$, is
continuous for all $v \in S-V_K^{\G,a}$.

Using translation under addition, we see that $\psi _v$ is also
uniformly continuous. Therefore, we can complete $\psi _v$ to
$\overline{\psi _v} :\overline{L^v} \rt K_{\tau (v)}$. Each
$\overline{\psi _v}$ is injective since any field homomorphism is
injective.

Now let $\overline {\ppr} : \mathcal{B} (G) \rt \mathcal{B} (G) $
be the embedding induced by the homomorphisms $\beta _v \circ
\overline{\psi _v} ^0: \G (\overline{L^v})\rt \G(K_{\tau (v)})$.
The map $\overline {\ppr}$ is clearly continuous on the
Furstenberg boundary, and since the Furstenberg boundary is
compact, $\oppr$ is uniformly continuous.

\ep

If $K$ is a number field then $\oppr$ is an automorphism. In
general though, it is not necessarily the case that a
self-embedding of a spherical building is an automorphism. Take
for example the spherical building for the standard flag complex
of $\mathbb{P}^k (\mathbb{F}_q ((t)))$ which is both isomorphic
to, and properly contains, the flag complex for $\mathbb{P}^k
(\mathbb{F}_q ((t^2)))$.

The surjectivity of $\oppr$ will be shown in Lemma 6.8 and must
wait until we can show that $\oppr$ extends $\pp$. Then we can use
the fact that $\pp$ has a dense image.

\bigskip \noindent \textbf{Extending the a.e. defined boundary function.} Our
goal is to show that $\pp$ is extended by $\oppr$.

Earlier we chose each global field $L^v$ to be large in an
algebraic sense with respect to each $K_v$. We can also assume
that each $L^v$ is topologically large with respect to each $K_v$
by choosing $L^v < K_v$ to be a dense subfield. Indeed, if $L^v$
is not dense we could replace $L^v$ with a finite extension that
is dense in $K_v$. This will ensure that
 $\mathcal{B} (G_R)$ carries some of the topological information
of $\mathcal{B} (G)$. In particular we have

\bigskip \noindent \textbf{Lemma 6.4} \emph{For any $(n-2)$-dimensional simplex
$f\subset \mathcal{B} (G_R)$, the set ${\mathcal{N}}_R(f) $ is
dense in  ${\mathcal{N}}(f) \se \mathcal{B} (G)$ under the
subspace topology of the Furstenberg topology.}

\pru Let $\Sigma _f \se \mathcal{B}(G_R)$ be an apartment
containing $f$, and suppose $c_f\se \Sigma _f$ is a chamber
containing $f$.

As in the proof of the previous lemma, there is a valuation $v \in
S-V_K^{\G,a}$ and an $L^v$-defined root subgroup $\mathbf{U}<\G$,
such that $\mathbf{U}(L^v)<G_R$ acts faithfully and transitively
on the set $\mathcal{N}_R(f)-\{\,c_f\,\}$. It also follows from
Proposition 5.6(i) of \cite{Ti}, that $\mathbf{U}(K_v)<G$ acts
faithfully and transitively on the set
$\mathcal{N}(f)-\{\,c_f\,\}$. Therefore, $\mathbf{U}(L^v)$ is
homeomorphic to $\mathcal{N}_R(f)-\{\,c_f\,\}$, and
$\mathbf{U}(K_v)$ is homeomorphic to $\mathcal{N}(f)-\{\,c_f\,\}$.

Since $L^v$ is dense in $K_v$, and because $\mathbf{U}$ is
isomorphic as an $L^v$-variety to affine space, we have that
$\mathbf{U}(L^v)$ is dense in $\mathbf{U}(K_v)$. Therefore, we
have the following series of dense inclusions
\begin{align*}\mathcal{N}_R(f)-\{\,c_f\,\} &
\se \mathcal{N}(f) -\{\,c_f\,\} \\
 & \se \mathcal{N}(f)
 \end{align*}

\ep

Let $F_R$ be the set of $(n-2)$-dimensional simplices in
$\mathcal{B} (G_R)$ and define $$D_R=\bigcup_{f \in
F_R}\mathcal{N}_U(f)$$

We use the topological properties of $\mathcal{B} (G_R)$, and of
$\oppr$, to deduce topological properties of $\pp |_{D_R}$ in the
following

\bigskip \noindent \textbf{Lemma 6.5} \emph{The function
$\pp |_{D_R}: D_R \rt \mathcal{B} (G)$ is Furstenberg continuous.}

\pru Let $\e > 0$ and a chamber $c_1 \subset D_R$ be given.

By Proposition 6.3, there is a $\delta _R
>0$ such that $$\widehat{d}\Big(\pp (w_1)\,,\, \pp (w_2)\Big)<\e /3$$ for all
chambers $w_1,w_2 \subset \mathcal{B} (G_R)$ with
$\widehat{d}(w_1,w_2) < \delta _R$.

Suppose $c_2 \subset D_R$ is a chamber with $\widehat{d}(c_1,c_2)<
\delta _R /3$. By Lemma 6.4 and Lemma 9.4, there are chambers
$c_i' \subset \mathcal{B} (G_R)$ that intersect $c_i$ in an
$(n-2)$-dimensional simplex, and such that $\widehat{d}(c_i,c_i')
< \delta _R /3$ and $\widehat{d}(\pp(c_i),\pp(c_i'))<\e /3$. Hence
$$\widehat{d}(\pp (c_1),\pp (c_2))\leq \widehat{d}(\pp (c_1'),\pp
(c_2'))+\Sigma_{i=1}^2 \widehat{d}(\pp (c_i),\pp (c_i'))<\e.$$

\ep

Since $ \pp |_{D_R}$ and $\oppr$ are continuous we have

\bigskip \noindent \textbf{Lemma 6.6} \emph{For any
simplex $q \subset D_R$, we have $\pp (q) = \oppr (q).$}

\pru Both $\pp |_{D_R}$ and $\oppr |_{D_R}$ are continuous so they
are uniquely determined by $\ppr$. Indeed, according to Lemma 6.4,
$\mathcal{B} (G_R)$ is Furstenberg dense in $D_R$.

\ep

In Section 8, a maximal $K_v$-split torus $\mathbf{A_v} <\G$ is
chosen for each $v \in S-V_K^{\G,a}$. The tori are used to supply
an ergodic theory argument that allows for the creation of the
boundary function $\pp :\pu \rt \mathcal{B} (G)$.

Let $\Sigma _A \subset \mathcal{B} (G)$ be the apartment
stabilized by the group $$\prod_{v\in S
-V_K^{\G,a}}\mathbf{A_v}(K_v)<G.$$ By conjugating $\ar$, we may
assume that $\Sigma _A$ is an apartment in $\Delta (G _R)$. Let
$W$ be the Weyl group with respect to $\Sigma _A$, and denote a
fixed chamber in $\Sigma_A$ by $a^+$. Let $a^-$ be the chamber in
$\Sigma_A$ opposite of $a^+$. For each $w \in W$ we let $P_w <G$
be the stabilizer of $w a^+$.

In Section 8, we will see that there exists a co-null subset
$\mathcal{U} \se G$ such that $\pu=\mathcal{U}a^+$. By Fubini's
theorem, we can conjugate $\ar$ such that $P_w \cap \mathcal{U}$
is co-null in $P_w$ for all $w\in W$.

Define $$U_\partial ^w = \{\,gw a ^- \in \widehat{X} \mid g\in P_w
\cap {\mathcal{U}} \, \}$$ and $${\mathcal{U}}^w = \{\,g\in
{\mathcal{U}} \mid gw a ^- \in {\mathcal{U}}_\partial ^w \,\}.$$
Note that $w a^-$ is opposite of $w a^+$, so we have that $P_w w
a^- $ is a full measure subset of $\widehat{X}$.  Since $P_w \cap
\mathcal{U}$ is co-null in $P_w$, it follows that $\pu ^w$
 is a full measure subset of $\widehat{X}$.
Hence, ${\mathcal{U}}^w \se G$ is co-null for all $w\in W$.
Consequently, $\cap_{w\in W}{\mathcal{U}}^w \se G$ is co-null.

We replace ${\mathcal{U}}$ with $\cap_{w\in W}{\mathcal{U}}^w$. As
a result, if $c \subset \pu$ is a chamber, then there is an
apartment $\Sigma _c$ which is completely contained in $\pu$, and
such that the chamber opposite from $c$ in $\Sigma _c$ is
contained in $\Sigma _A$. For any chamber $c\subset \pu$, we let
$$\delta_A(c)=\min_{\Sigma _c} \{\, d_{\Sigma _c} (c, \Sigma
_A) \,\},$$ where the $\min$ is taken over all $\Sigma _c$ as
above with respect to the Tits metric $d_{\Sigma _c}$ on $\Sigma
_c$.

We can now improve upon Lemma 6.6.

\bigskip \noindent \textbf{Lemma 6.7} \emph{For any
simplex $q \subset \pu$, we have $\pp (q) = \oppr (q).$}

\pru For a chamber $c \subset \pu$, we prove that $\pp (c) = \ppr
(c)$ by induction on $\delta_{A}(c)$.

If $\delta_{A}(c)\leq 1$, then the result follows from the
previous lemma. Now suppose the result is true for any chamber $f
\subset \pu$ with $\delta_{A}(f)\leq k-1$, and let $c\subset \pu$
be a chamber with $\delta_{A}(c)=k$.

Let $\Sigma _c \subset \pu$ be an apartment containing $c$, and
such that the chamber in $\Sigma _c$ opposite of $c$ is contained
in $\Sigma _{A}$. Choose a chamber $f\subset \Sigma _c$ such that
$d_{\Sigma _c} (c, f)=1$ and $\delta_{A}(f)<k$. If $f^{op}$ is the
chamber in $\Sigma _c$ opposite of $f$, then $\de_{A}(f^{op})\leq
1$. By our induction hypothesis, $\pp (f) = \ppr (f)$ and $\pp
(f^{op}) = \ppr (f^{op})$.

It will be shown in Lemma 8.8 that $\pp$ preserves apartments.
Therefore, $\pp (\Sigma _c)$ is an apartment. In fact, $\pp
(\Sigma _c)$ is the unique apartment containing $\ppr (f)$ and
$\ppr (f^{op})$. Note that $\ppr (\Sigma _c) $ is also the unique
apartment containing $\ppr (f)$ and $\ppr (f^{op})$.

We conclude our proof by observing that both $\pp (c)$ and $\ppr
(c)$ must be the unique chamber in $\pp (\Sigma _c) = \ppr (\Sigma
_c)$ that contains $\pp (c\cap f) =\ppr (c \cap f)$, but not $\pp
(f)=\ppr (f)$.

\ep

\bigskip \noindent \textbf{The extension is an automorphism.}
 Now that we have shown that $\oppr$ extends $\pp$, we have to
prove that $\oppr$ is surjective, and hence an automorphism of
$\mathcal{B} (G)$. Then it follows that $\oppr$ corresponds to an
automorphism of $G$, or alternatively, an isometry of $X$.

\bigskip \noindent \textbf{Lemma 6.8} \emph{
The map $\overline{\ppr}$ is an automorphism of} $\mathcal{B}
(G)$.

\pru Let $\p^*$ be a coarse inverse for $\p$, and define $\pu*$
and $\overline{\pp _R^*}$ analogously to $\pu$ and $\oppr$.

Let $\Sigma \in \pu^*$, and let $F \se X$ be the flat
corresponding to $\Sigma$. Note that $\p \circ \p^*$ preserves the
portion of $F$ that lies near an orbit of $\ar$ in $X$ (see
Section 8). Since $F$ is the only flat in $X$ that is a finite
Hausdorff distance from itself, it follows that
$$\oppr \circ \overline{\pp _R^*}
(\Sigma) = \Sigma.$$ Hence,
$$\pu ^* \se \oppr \big(\mathcal{B}(G)\big).$$

Note that the map $\oppr$ either has a closed null image or is
surjective since $K_{\tau(v)} $ is a $\overline{\psi
_v}(\overline{L^v})$-vector space. The lemma follows since $\pu
^*$ is co-null in $\widehat{X}$.

\ep

\bigskip \noindent \textbf{Automorphisms that correspond to quasi-isometries.}
Let $\Hd$ denote the Hausdorff distance between closed subsets of
$G$. We define the group  $$\aut_{\Hd}(G ;\ar)=\{\, \varphi \in
\aut(G) \mid \Hd(\varphi (\ar)\,,\,\ar)< \infty \,\}.$$

Using Lemma 6.7 and Lemma 8.3(vii), Eskin's proof that the
automorphism $\oppr \in \aut(G) \cong \text{Isom}(X)$ corresponds
to an isometry of $X$ that is a finite distance from $\p$
(\cite{Es} Step 7) can be applied to show

\bigskip \noindent \textbf{Proposition 6.9} \emph{There is an
isomorphism of topological groups }$$\mathcal{QI}(\ar) \cong \aut
_{\Hd}(G ; \ar)$$

\bigskip

The proof proceeds by identifying points in $X$ as intersections
of flats in $X$. Flats are parameterized by apartments in
$\mathcal{B}(G)$, so $\oppr$ completely determines where points in
$X$ are mapped to under the corresponding isometry of $X$. Any
point in a $\ar$-orbit is a bounded distance from the intersection
of flats whose boundaries are in $\pu$. Therefore, $\p$ maps
points in a $\ar$-orbit to within a bounded distance of their
images under the isometry corresponding to $\oppr$.

Eskin's proof makes no mention of the topological nature of this
isomorphism, but it clearly follows. The fact that the isomorphism
is topological is more interesting in the $S$-arithmetic setting
since merely the abstract group type of the quasi-isometry group
of an arithmetic lattice in a real semisimple Lie group determines
the lattice up to commensurability.

\bigskip \noindent \Large \textbf{7. Automorphisms coarsely preserving lattices} \normalsize \bigskip

We want to determine the group $\aut_{\Hd}(G;\ar)$ and complete
our proof of Theorem 1.1.

\bigskip \noindent \textbf{The case of anisotropic groups.}
Notice that if $\G$ is $K$-anisotropic, then $\aut_{\Hd}(G;\ar)$
is isomorphic to $\aut (G)$. Indeed, $\ar$ is a cocompact lattice
in $G$ so $\Hd(G \,,\,\ar) < \infty$. Thus, our proof of Theorem
1.1(iii) is complete (assuming the results from Sections 8 and 9).

\bigskip \noindent \textbf{The function field case for isotropic groups.}
The proof of Theorem 1.1(ii) concludes with Lemma 7.1 below. We
include the proof here to group it with similar results, but its
proof uses notation and concepts defined in Section 8. The reader
may want to return to the proof of this small fact after having
read what will follow.

\bigskip \noindent \textbf{Lemma 7.1} \emph{If $\G$ is $K$-isotropic,
 then the group $\aut_{\Hd}(G;\ar)$ is a measure
zero subgroup of $\aut(G)$.}

\pru For a given element of $\aut_{\Hd}(G;\ar)$, we let $g:X\rt X$
be the corresponding isometry. We choose a neighborhood
$N(\Gamma)^g \se X$ of the set $\neu$ from Lemma 8.3, such that
$\neu \se g(\neu ^g)$.

Define $\vol_F$ to be Lebesgue measure on $F$, and let $\e$ be as
in Lemma 8.3. There is a Weyl chamber $\mathfrak{C} \se X$ such
that for any $g \in \aut_{\Hd}(G;\ar)$, for any flat $F \se X$
that contains $\mathfrak{C}$ up to Hausdorff equivalence, and for
any point $x\in F$, we have
$$ \lim_{r \to \infty} \frac{\vol_{F}\Big(\big[F \cap
N\big(\Gamma \big)^g \big] \cap
B_x(r)\Big)}{\vol_{F}\Big(B_x(r)\Big)} < 1-\varepsilon .$$

Let $F'\se X$ be a flat containing $g(\mathfrak{C})$ up to
Hausdorff equivalence. Then, by replacing $F$ with $g^{-1}(F')$ in
the preceding inequality, it follows that for any point $y \in
F'$:
$$ \lim_{r \to \infty} \frac{\vol_{F'}\Big(\big[F' \cap
N\big(\Gamma \big) \big] \cap
B_y(r)\Big)}{\vol_{F'}\Big(B_y(r)\Big)} < 1-\varepsilon . $$
Therefore, $F' \notin \mathfrak{U}$. Hence, if $c \subset
\mathcal{B}(G)$ is the chamber representing the equivalence class
of $\mathfrak{C}$, then $\aut_{\Hd}(G;\ar)\cdot c \se
\mathcal{B}(G) - \pu$.

The lemma follows from Fubini's theorem since $\pu $ is co-null in
$\widehat{X}$.

\ep

\bigskip \noindent \textbf{The number field case for isotropic groups.}
The proof of the following proposition was indicated to me by
Nimish Shah, and it completes the proof of Theorem 1.1(i).

\bigskip \noindent \textbf{Proposition 7.2} \emph{If $K$ is a
number field and $\G$ is $K$-isotropic, then} $$\aut_{\Hd}\big(G;
\ar \big) = \comm\big(\ar\big).$$

\pru Let $\varphi \in \aut_{\Hd}(G,\ar)$. We have to show that
$\varphi \in \comm(\ar)$.

To simplify notation we let $\Lambda = \ar$ and $\Lambda ^\varphi
= \varphi (\ar)$. By replacing $\Lambda$ with a finite index
subgroup, we can assume that $\Lambda$ and $\Lambda ^\varphi $ are
contained in the group $$G ^+=\prod _{v \in S} \G(K_v)^+.$$

By Ratner's theorem on unipotent flows (\cite{Ra} Theorem 6.4),
the orbit of the point $(\Lambda , \Lambda ^\varphi )$ in $ G^+ /
\Lambda \times \hp / \Lambda ^\varphi $ under the diagonal action
of $\hp$ is homogeneous. If we denote the diagonal embedding of
$G^+$ into $\hp \times \hp$ by $\Delta G^+$, then the previous
sentence says that
$$\overline{\Delta \hp (\Lambda , \Lambda ^\varphi)}=L(\Lambda ,\Lambda ^\varphi),$$ where $L$
is a closed subgroup of $\hp \times \hp$ which contains $\Delta
G^+$.

We claim that either $L= \Delta \hp $ or there is some $v\in S$
such that $$1 \times {\bold{G}}(K_{v}) ^+ \leq (1 \times \hp) \cap
L.$$ Indeed, if $\Delta \hp <L$, then there are group elements
$g_1,g_2 \in \hp$ such that $(g_1,g_2) \in L$ and $g_1 \ne g_2$.
Hence, there is some $g \in \hp$ with $g \ne 1$ and $(1,g)\in L$.
That is to say, $(1 \times \hp)\cap L$ is nontrivial. Note that if
$(1,h) \in L$, then for any $g \in \hp$, we have $(1,
ghg^{-1})=(g,g)(1,h)(g^{-1},g^{-1})\in L$ since $\Delta \hp <L$.
Thus, $(1 \times \hp) \cap L$ is a normal subgroup of $1 \times
\hp$. Now a theorem of Tits' (\cite{Ti 1} Main Theorem) tells us
that that each group $\G(K_v)^+$ is simple since $\G$ has a
trivial center. Therefore, $1 \times {\bold{G}}(K_{v}) ^+ \leq L$
for some $v \in S$, since $(1 \times \hp) \cap L$ is a nontrivial
normal subgroup of $1 \times \hp$. Thus, our claim is proved.

If it is the case that $1 \times {\bold{G}}(K_{v}) ^+ \leq L$,
then $ \overline{\Delta \hp (\Lambda , \Lambda ^\varphi)}$
contains $\{ \Lambda \} \times \hp / \Lambda ^\varphi $, as
$\Lambda ^\varphi$ is irreducible. Hence, for any $g \in \hp$,
there is a sequence $\{ g_k \} \se \hp$ such that $\Delta g_k
(\Lambda ,\Lambda ^\varphi ) \to (\Lambda, g \Lambda ^\varphi ).$
Since $g _k \Lambda \to \Lambda$, it follows that there are
sequences $\{h_k\} \se \hp$ and $\{\lb _k\} \in \Lambda$, such
that $g_k = h_k \lb _k$ and $h_k \to 1$. Therefore, $h_k^{-1} g_k
\in \Lambda$ and $h_k^{-1}g_k \Lambda ^\varphi \to g \Lambda
^\varphi$ which proves that $\overline{\Lambda \Lambda ^\varphi
}=\hp / \Lambda ^\varphi$. Note that our assumption that $\varphi
\in \aut_{\Hd}(G;\Lambda )$ implies that $\Lambda \Lambda
^\varphi$ is bounded. Thus, this case is precluded.

We are left to consider the case when $\Delta \hp =L$. We will
show that $\Lambda \Lambda ^\varphi \se \hp / \Lambda ^\varphi$ is
a closed set. To this end, suppose there is a sequence $\{\lb _k\}
\se \Lambda$ and a group element $g\in \hp$ with $\lb _k \Lambda
^\varphi \to g \Lambda ^\varphi$. Then $\Delta \lb _k (\Lambda ,
\Lambda ^\varphi ) \to (\Lambda , g \Lambda ^\varphi)$. Since
$\Delta \hp (\Lambda , \Lambda ^\varphi)$ is closed, $(\Lambda , g
\Lambda ^\varphi )= \Delta h (\Lambda , \Lambda ^\varphi )$ for
some $h \in \hp$. Therefore, $g \Lambda ^\varphi = h \Lambda
^\varphi$. Since $h \Lambda = \Lambda $, we have $h \in \Lambda$
which shows that $\Lambda \Lambda ^\varphi $ is closed.

Since $\Lambda \Lambda ^\varphi$ is bounded, it must be compact
which would require it to be finite or perfect. As perfect sets
are known to be uncountable, $\Lambda \Lambda ^\varphi$ is finite.
That is $\varphi \in \comm(\ar)$ as desired.

\ep

Assuming the material from Sections 8 and 9, the proof of Theorem
1.1 is complete. It is the absence of the counterpart to
Proposition 7.2 for function fields that leads to the discrepancy
between (i) and (ii) of Theorem 1.1 and Corollary 1.2.

\bigskip \noindent \textbf{The commensurator group.} We close
 this section with a lemma that provides a
concrete description of $\comm(\ar)$.

\bigskip \noindent \textbf{Lemma 7.3} \emph{The group $\comm (\ar)$
is an extension of $\Aut(\G)(K)$ by $\aut(K)_{\G,S}$. If $\G$ is
$K$-split and defined over a subfield of $K$ that is fixed
pointwise by $\aut(K)_{\G,S}$, then}
$$\comm(\ar) \cong \Big( \G (K) \rtimes  \out (\G)(K)
\Big) \rtimes \aut(K)_{\G,S}.$$

\pru Recall that $\ar$ is embedded diagonally in $G$ with respect
to the simple factors of $G$. Hence, any group element in $\comm
(\ar) \cap G$ would have to take a finite index diagonal subgroup
of $\ar$ into the diagonal of $G$. It follows from the Borel
density theorem that any finite index subgroup of $\ar$ is a
Zariski dense subset in each simple factor of $G$. Therefore,
$\comm (\ar) \cap G$ is also contained in the diagonal of $G$.

We have shown that, as an abstract group, $\comm (\ar) \cap G$ is
a subgroup of the group $L$ of inner automorphisms of $\G (K_v)$
which commensurate $\ar < \G(K_v)$; the choice of $v \in
S-V_K^{\G,a}$ is arbitrary.

Borel's well known determination of inner commensurators for
arithmetic groups (\cite{Bo} Theorem 2) essentially contains a
proof that $L=\G (K) < \G(K_v)$. Therefore, $\comm (\ar) \cap G$
is the diagonal subgroup $ \Delta \G (K) <G$.

If $\varphi \in \comm(\ar )$, then $\ar$ and $\varphi(\ar)$ are
commensurable. Hence, an inner automorphism of $G$ commensurates
$\ar$ if and only if it commensurates $\varphi(\ar)$. Therefore,
$\varphi ( \Delta \G (K))= \Delta \G(K)$.

Conversely, suppose $\varphi $ is an automorphism of $G$ with
$\varphi (\Delta \G (K))= \Delta \G (K)$. Then $\varphi (\ar )$ is
a lattice contained in $\Delta \G (K)$, so $\varphi (\ar)$ is
commensurable to $\ar$ by the proof of the Margulis-Venkataramana
arithmeticity theorem (see \cite{Mar} pages 307-311). Therefore,
$\varphi \in \comm(\ar)$.

Hence, finding $\comm(\ar)$ amounts to finding the subgroup of
$\aut(G)$ that stabilizes $\Delta \G (K)$. This is what we shall
do.

Suppose $\psi \in \aut(G)$ and that $\psi (\Delta\G(K))=\Delta
\G(K)$. By Theorem (A) of \cite{Bo-T}, $\psi \in \aut(G)$ can be
uniquely written in the form
$$\prod _{v \in S-V_K^{\G,a}} \beta _v \circ \alpha _v ^\circ $$ for some permutation $\tau$ of
$S-V_K^{\G,a}$, a collection of field isomorphisms $\alpha _v :K_v
\rt K _{\tau (v)}$, and a collection $\beta _v :
\sideset{^{\alpha_v}}{}{\G} \rt \G$ of $K _{\tau
(v)}$-isomorphisms of algebraic groups. Since $\psi$ is a
homeomorphism, each field isomorphism $\alpha_v$ is a
homeomorphism as well.

Since $\Delta \G (K)$ is stabilized by $\psi$, $$\beta_v \circ
\alpha _v ^\circ |_{\G (K)}=\beta_w \circ \alpha _w ^\circ |_{\G
(K)}$$ for all $v,w \in S-V_K^{\G,a}$. Again by Theorem (A) of
\cite{Bo-T}, there exists a unique $\sigma \in \aut (K)$ and a
unique $K$-isomorphism of algebraic groups $\delta :
\sideset{^{\sigma}}{}{\G} \rt \G$, such that $\delta \circ \sigma
^0$ is extended by all $\beta _v \circ \alpha _v ^\circ$.

Because each $ \alpha _v $ is a homeomorphism, $\sigma :K \rt K$
is a homeomorphism between $K$ with the $v$-topology and $K$ with
the $\tau(v)$-topology. Therefore, $\tau (v) = \sigma \cdot v$ for
all $v \in S-V_K^{\G,a}$. That is, $\sigma \in \aut(K)_{\G,S}$.

We have identified an inclusion of $\comm (\ar)$ into the group of
pairs $(\delta , \sigma )$, where $\sigma \in \aut(K)_{\G,S}$ and
$\delta : \sideset{^{\sigma}}{}{\G} \rt \G$ is a $K$-isomorphism.
To see that the inclusion is an isomorphism, let $(\delta , \sigma
)$ be a given pair as above. For any $v \in S-V_K^{\G,a}$, let
$\sigma _v :K \rt K$ be defined by $\sigma_v(x)=\sigma(x)$. We
assume that the domain of $\sigma _v$ has the $v$-topology and
that the image of $\sigma _v$ has the $\sigma \cdot v$-topology.
Hence, $\sigma _v$ is continuous, and it may be completed
topologically to obtain an isomorphism $\overline{\sigma _v} :K_v
\rt K_{\sigma \cdot v}$. Then we define a homomorphism $\G (K_v)
\rt \G (K_{\sigma \cdot v})$ by $\delta \circ \overline{\sigma _v}
^\circ$. The product map
$$\prod _{v\in S-V_K^{\G,a}} \delta \circ \overline{\sigma _v}
^\circ$$ is then an automorphism of $G$ that stabilizes $\Delta \G
(K)$. Hence, the group of pairs $(\de , \sigma)$ as above is
isomorphic to $\comm (\ar)$.

Notice that the group operation on $\comm (\ar)$ is given by $(\de
, \sigma)(\de ' ,\sigma ')=(\de \circ \sideset{^\sigma}{}{\del}',
\sigma \sigma ')$, where $\sideset{^\sigma}{}{\del}' :
\sideset{^{\sigma  \sigma '}}{}{\G} \rt \sideset{^\sigma}{}{\G}$
is the $K$-isomorphism obtained by applying $\sigma $ to the
coefficients of the polynomials defining $\delta '$. This is the
group structure of an extension:
$$1 \rt \mathbf{Aut(G)}(K) \rt \comm(\ar) \rt \aut (K)_{\G,S}
\rt 1. $$

The above extension splits if $\G$ is defined over a subfield of
$K$ that is fixed pointwise by $\aut(K)_{\G,S}$. Indeed, if $\G$
is defined over such a field, then for any $\sigma \in \aut(K)$ we
have $\sideset{^\sigma}{}{\G}=\G$. It follows that if $id_{\G}:\G
\rt \G$ is the identity map, then the pairs $(id_{\G},\sigma)$
exist in $\comm(\ar)$. Hence, the extension splits.

For the statement that $\G$ being $K$-split implies
$$\mathbf{Aut(G)}(K) \cong \G (K) \rtimes \mathbf{Out (G)}(K),$$
see, for example, the discussion in 5.7.2 of \cite{Ti}. (Recall
that we identify $\mathbf{Out (G)}$ with the automorphism group of
the Dynkin diagram of $\G$.)

\ep

\noindent \Large \textbf{8. Constructing a boundary function
defined a.e.} \normalsize \bigskip

Sections 6 and 7 show the conclusion of the proof for Theorem 1.1
once the boundary function $\pp : U_\partial \rt \mathcal{B} (G)$
is created. In section 8, we outline the construction of $\pp$. We
will refer to \cite{Es} for most of the details of the
construction.

\bigskip \noindent \textbf{Replacing the word metric.} Let $\mathbf{\widetilde{G}}$ be the algebraic simply connected
cover of $\G$. We define $$H = \prod _{v \in S-V_K^{\G,a}}
\mathbf{\widetilde{G}}(K_v)$$ and $$\Gamma =
\mathbf{\widetilde{G}}(\mathcal{O_S}).$$ Note that $\Gamma $ and
$\ar$ are commensurable up to finite kernels (see e.g. \cite{Mar}
I.3.1.1.iv).

Let $\mathfrak{K}$ be a maximal compact subgroup of $\tg$, and let
$\varepsilon ' >0$ be given. Let $\mu$ be the probability measure
on $\Gamma \backslash H$ which is derived from Haar measure on
$H$. We choose a compact set $D \se \Gamma \backslash H$ which
contains the coset $\Gamma$, and such that $\mu (D) \geq
1-\varepsilon '$.

We denote by $N(\ga)^\circ \subseteq \tg /\mathfrak{K}$ the set of
all cosets with a representative in $H$ that maps into $D$ under
the quotient map $H \rt \Gamma \backslash H$. In symbols,
$$N(\Gamma )^0 = \{\, h\mathfrak{K} \in \tg /\mathfrak{K} \mid \tga h \in D \,\}. $$

Since $\mathfrak{K}$ is the isotropy group of a point in $X$, we
can identify $H/\mathfrak{K}$ as a subset of $X$. For each
$h\mathfrak{K} \in \tg /\mathfrak{K} $, we let $P(h\mathfrak{K})$
be the set of points in $X$ that are at least as close to
$h\mathfrak{K} \in X$ as to any other point of $H/\mathfrak{K} \se
X$. Precisely:
$$P(h\mathfrak{K})=\{\, x\in X \mid d (x,h\mathfrak{K})\leq
d(x,g\mathfrak{K}) \text{ for all } g\in \tg \,\}.$$

Let $$N(\ga)=\bigcup_{h\mathfrak{K} \in N(\ga)^\circ}
P(h\mathfrak{K}).$$ Notice that $N(\Gamma ) \se X$ contains the
orbit $\Gamma \mathfrak{K}$. Since $\ga \backslash N(\ga)^\circ =D
$ and $P(\mathfrak{K})$ are compact, $\ga \backslash N(\ga) =\ga
\backslash [N(\ga )^\circ P(\mathfrak{K})]$  is compact. Thus,
$\Gamma$ is quasi-isometric to $\neu \se X$ with the path metric.

The geometry of $\neu \se X$ with the path metric is more
convenient to work with than the word metric on $\Gamma$. More
convenient still, would be working with the geometry of $\neu$
under the restricted metric from $X$.

In general, a lattice is not quasi-isometric to its orbit with the
restricted metric, but with our standing assumption that
$\mathbf{\widetilde{G}}$ is placewise not rank one, we can apply
the theorem below from \cite{L-M-R}

\bigskip \noindent \textbf{Theorem 8.1
(Lubotzky-Mozes-Raghunathan)} \emph{The word metric on $\Gamma$ is
quasi-isometric to $\neu \se X$ with the restricted metric.}

\bigskip

Using Theorem 8.1, the fact that $\Gamma$ and $\ar$ are
commensurable up to finite kernels, and the fact that the
inclusion of $\neu$ with the restricted metric into $X$ is
isometric, we can realize a given quasi-isometry $$\p : \ar \rt
\ar$$ by a quasi isometric embedding $$\neu \rt X.$$ The resulting
embedding is a finite distance in the sup norm from $\p $, so we
will also denote it by $\p$. We will assume that $$\p : \neu \rt
X$$ is a $(\kappa ,C)$ quasi-isometric embedding.

\bigskip \noindent \textbf{Ergodic actions of abelian groups.}
For each $v \in S-V_K^{\G,a}$, let ${\bf A_v}$ be a maximal
$K_v$-split torus in $\mathbf{\widetilde{G}}$. We define the group
$$A=\prod_{v \in S-V_K^{\G,a}} {\bf A_v}(K_v)<H.$$ We denote the flat
 corresponding to $A$ by $\ap \se X$. We may assume that
$\mathfrak{K} \in \ap$.

We introduce a pseudometric $d_{A}$ on $A$ by setting
$d_A(a_{1},a_{2})$ to be equal to
$d(a_{1}\mathfrak{K},a_{2}\mathfrak{K})$ for $a_1 \mathfrak{K},
a_2 \mathfrak{K} \in X$.

There is also a Haar measure on $A$ which we denote by $da$. We
denote Lebesgue measure on $\ap$ by $\vol _\ap$. Then, after a
normalization, we have for any measurable set $Y \subseteq A$:
$$da\Big(A\cap \Big(\bigcup _{a\in Y}a\mathfrak{K}a^{-1}\Big)\Big)=
\vol_{\ap}\Big(\ap \cap \Big(\bigcup _{a\in Y}
aP(\mathfrak{K})\Big)\Big).$$

The Birkhoff ergodic theorem is usually stated for an ergodic
action of $\mathbb{Z}$. However, a careful reading of the proof of
the Birkhoff ergodic theorem shows that it applies to ergodic
actions of our pseudometric group $A$ as well (see e.g. \cite{Bl}
Theorem 3.2). That is, if we let $B_1^A(r) \se A$ be the ball of
radius $r$ centered at the identity element of $A$, then we have
the following

\bigskip \noindent \textbf{Proposition 8.2 (Birkhoff ergodic
 theorem)} \emph{If $Y$ is a finite
volume right ergodic $A$-space and $f\in L^1 (Y)$, then for a.e.
$y\in Y$:}
$$\lim_{r \to \infty} \frac{1}{da(B_1^A(r))}\int _{B_1^A(r)}f(ya)da
=\int_{Y} f .$$

\bigskip

Prasad's proof of the strong approximation theorem for simply
connected semisimple Lie groups contains a proof of the ergodicity
of the $A$-action on $\Gamma \backslash H$ (see \cite{Pr 2} Lemma
2.9). Hence, we can apply the Birkhoff ergodic theorem to the
action of $A$ on $\Gamma \backslash H$.

\bigskip \noindent \textbf{Generic flats have most of their
volume near $\Gamma$.} Following Eskin, we are now prepared to
show that a generic flat in $X$ has most of its volume contained
in $\neu \se X$.

For any group element $h\in H$, define  $\vol_{h \ap}$ to be
Lebesgue measure on the flat $h \ap \se X$. That is, for any
measurable set $Y\se h \ap$, we let $$\vol _{h\ap}(Y)=\vol _\ap
(h^{-1} Y).$$ Thus, the measure $\vol _{h\ap}$ is compatible $da$
in a natural way.

We denote by $B_x^{h\ap}(r)\se h\ap$ the metric ball centered at
the point $x\in h\ap$ with radius $r>0$. Denote the characteristic
functions of $\neu \se X$ and $D \se \Gamma \backslash \neu$ by
$\chi _{\neu}$ and $\chi _D$ respectively.

By Proposition 8.2, we have that for $\mu$ a.e. $\tga h \in
\gmga$:

\begin{align*} \lim_{r \to \infty}& \frac{1}{\vol_{h \ap }(B ^{h \ap}_{h\mathfrak{K}}(r))}
 \int _{B^{h \ap}_{h\mathfrak{K}}(r)} \chi _{\neu} \vol_{h \ap} \\
 & =  \lim_{r \to \infty} \frac{1}{da(B^{A}_{1}(r))}
\int _{B^{A}_{1}(r)} \chi _{\neu}(ha\mathfrak{K}) da     \\
  & \geq  \lim_{r \to \infty} \frac{1}{da(B^{A}_{1}(r))}
\int _{B^{A}_{1}(r)} \chi _{D}(\tga ha) da \\
 &  =  \int _{\gmga} \chi _{D} \\
 & = \mu (D) \\
 & \geq 1-\varepsilon '.
\end{align*}
The inequality shows that for a.e. $\tga h \in \gmga$, any $\gamma
\in \Gamma $, and any point $x\in \gamma h \ap$:
$$ \lim_{r \to \infty} \frac{\vol_{\gamma h \ap}\Big(\big[\gamma h \ap \cap
N\big(\Gamma \big)\big] \cap B^{\gamma h
\ap}_x(r)\Big)}{\vol_{\gamma h \ap}\Big(B^{\gamma h
\ap}_x(r)\Big)} \geq 1-\varepsilon '.$$ Hence, the generic flat
has much of its volume contained in $\neu$.

The above argument is the basic idea behind Lemma 8.3 below.
Refining the argument will yield more precise information about
how much of a generic flat is contained in $\neu$. Then we will be
in a position to apply the quasiflats with holes theorem from
\cite{W} to begin constructing a map on $\mathcal{B}(G)$.

\bigskip \noindent \textbf{More on the position of a generic flat with
 respect to $\Gamma$.} Let $h\in H$. For a set $W \se X$ contained
  in the flat $h \ap$, we let
$$W _{(\e , \rho)} = \{\, x\in W \mid B_y^{h\ap} \big(\e d(x,y)\big)\cap W
\neq \emptyset \mbox{ for all } y\in h \ap - B_x ^{h\ap}(\rho)
\,\}.$$ Hence, $W _{(\e , \rho)}$ is the set of all points $x\in
W$ which can serve as an observation point from which all points
in $h \ap$ (that are a sufficient distance from $x$) have a
distance from $W$ that is proportional to their distance from $x$.

We denote the metric $r$-neighborhood of a set $Y \se X$ by $\nd
_r (Y)$. We denote the Hausdorff distance between two sets $P,Q
\se X$ by $\hd (P,Q)$.

Recall the definition of a \emph{wall} $L \se X$ as a codimension
$1$ affine subspace of a flat, that is contained in at least two
distinct flats.

Lemma 8.3 below is an amalgam of Lemmas 2.2, 3.2, and 5.2 from
\cite{Es}. We omit the proof of the lemma as it is nearly
identical to those in \cite{Es}. We note that the proof follows
the principle shown above using the Birkhoff ergodic theorem.

 We will assume throughout that $\varepsilon >0$ is a sufficiently
 small number
 depending on $\kappa$ and $X$.

\bigskip \noindent \textbf{Lemma 8.3} \emph{There are constants
 $\rho >0$, and $\rho '
>0$ depending on $\e$ and $X$; constants  $\lambda_0>1$, $\lambda_1>1$,
  $N'>0$, $m> 0$, and $1>b>0$ depending on $X$; and a
$\Gamma$-invariant co-null set ${\mathcal{U}}\subseteq \tg$ such
that for any $h \in {\mathcal{U}}$ there are sets  $$\Omega _{h
\ap} ^* \subseteq \Omega _{h \ap}  \subseteq \Omega _{h \ap} '
\subseteq h \ap \cap \neu$$ which satisfy the following
properties:}

\emph{(i) For any point $x \in h \ap$: $$ \lim_{r \to \infty}
\frac{\vol_{ h \ap}\Big( \big[h \ap \cap N\big(\Gamma \big)\big]
\cap B^{ h \ap}_x(r)\Big)}{\vol_{ h \ap}\Big(B^{ h \ap}_x(r)\Big)}
\geq 1-\varepsilon/4.
 $$}

\emph{(ii) $\Omega _{ h \ap}'\subseteq \big( h \ap \cap
\neu\big)_{(\varepsilon, \rho)}$ and for any point $x \in  h \ap$:
$$ \lim_{r \to \infty} \frac{\vol_{ h \ap}\Big(\Omega _{ h \ap}'
\cap B^{ h \ap}_x(r)\Big)}{\vol_{ h \ap}\Big(B^{ h \ap}_x(r)\Big)}
\geq 1-\varepsilon/2.
 $$}

\emph{(iii) $\Omega _{ h \ap}\subseteq (\Omega _{ h
\ap}')_{(\varepsilon, \rho ')}$ and for any point $x \in  h \ap$:
$$ \lim_{r \to \infty} \frac{\vol_{ h \ap}\Big(\Omega _{ h \ap}
\cap B^{ h \ap}_x(r)\Big)}{\vol_{ h \ap}\Big(B^{ h \ap}_x(r)\Big)}
\geq 1-\varepsilon/2.
 $$}

\emph{(iv)  For any point $x \in  h \ap$: $$ \lim_{r \to \infty}
\frac{\vol_{ h \ap}\Big(\Omega _{ h \ap}^* \cap B^{ h
\ap}_x(r)\Big)}{\vol_{ h \ap}\Big(B^{ h \ap}_x(r)\Big)} \geq 1-m
\varepsilon.
 $$ }

\emph{(v) If $y \in \Omega _{ h \ap}^*$ and $L \subset  h \ap$ is
a wall with $d (y,L)<N'$ then there is a group element $h' \in
{\mathcal{U}}$ such that $$ h' \ap \cap  h \ap \subseteq
\Nd_{2N'}(L)$$ and
$$\Hd\Big( h \ap \cap \Nd_r( h' \ap)\;,\; L\Big) \leq \lambda _1 r$$ for any
$r>\lambda _0$.}

\emph{(vi) For any wall $L \subset \ap$ and any point $x\in h
\ap$:
$$ \lim_{r \to \infty} \frac{\vol_{ h \ap}\Big(\Omega _{ h\ap}^*
\cap B^{ h \ap}_x(r)\cap  h L \Big)}{\vol_{ h \ap}\Big(B^{ h
\ap}_x(r)\Big)} \geq b.
 $$ }

\emph{(vii) There is a $\Gamma$ invariant set $E \subseteq
{\mathcal{U}}$ such that $\mu (\Gamma \backslash E)>1-\varepsilon
/2$, and $h\mathfrak{K} \in \Omega _{h \ap}$ for any $h\in E$. }

\bigskip

\bigskip \noindent \textbf{Remarks.} There are some differences
 in this lemma with Lemmas 2.2, 3.2, and 5.2 in
\cite{Es}. In particular, the transverse flats in part (v) do not
necessarily intersect in a wall for the general space $X$, as can
be arranged if $X$ is a symmetric space. Take for example a
regular trivalent tree which is the Euclidean building for
${\bf{SL_2}}(\mathbb{Q}_2)$. The walls in this example are
vertices; the flats are lines, and there is no pair of lines which
intersect in a single point.

Also, the constant $b$ in part (vi) is shown in \cite{Es} to be
nearly one. This discrepancy is essentially due to the fact that
if $X$ is a Euclidean building, then the orbit of
$P(\mathfrak{K})$ under the action of the $\mathfrak{p}$-adic
group that stabilizes a wall containing $\mathfrak{K}$ may not
contain all of $L$. Take for example the building for $\sl3
(\mathbb{Q}_p)$. However, Eskin's proof only uses that the
constant is greater than $0$, and that is all we shall need as
well.

\bigskip \noindent \textbf{A collection of useful flats.} Lemma
8.3 provides us with a collection of flats in $X$ that have most
of their volume, and a substantial portion of the volume of their
walls, contained in $N(\Gamma)$. We denote this collection of
flats by $\mathfrak{U}$. That is,
$${\mathfrak{U}}=\{ h \ap | h \in {\mathcal{U}} \}.$$

Since any flat $F\in  \mathfrak{U}$ has most of its volume
contained in $\neu$, we can restrict $\p : \neu \rt X$ to $F \cap
\neu$ and begin to analyze the image using Theorem 1.2 of
\cite{W}. We state this theorem as

\bigskip \noindent \textbf{Theorem 8.4 (Quasiflats with holes)}
\emph{Let $\varphi :\Omega \rt X$ be a $(\kappa , C)$
quasi-isometric embedding of a set $\Omega \se \mathbb{E}^n$.
There are constants $M=M(\kappa ,X)$ and $\delta _0 = \delta _0
(\kappa ,X)$ such that if $\delta < \delta_0$, then there exists
flats $F_1,F_2,...,F_M \se X$ such that
$$\varphi \big(\Omega _{(\delta , R)}\big)\se \nd _N \Big(\bigcup_{i=1}^M F_i\Big),$$
where $N=N(\kappa ,C ,R ,X)$.}

\bigskip

Theorem 8.4, and the fact that a generic flat $F \se X$ is
contained in $\mathfrak{U}$, positions us to begin constructing
the function $\pp : \pu \rt \mathcal{B}(G)$ where the set $\pu \se
\mathcal{B}(G)$ has full measure in $\widehat{X}$.

\bigskip \noindent \textbf{Weyl chambers are mapped to Weyl chambers.}
For points $x,z,w \in X$ and a number $\rho \geq 0$, we let
$$D_x ( \rho ; z,w ) = \max \{ \rho, d(x, z), d(x,w)\}.$$ Define a function $\p : X \rightarrow Y$ to be a
$(\ka ,\rho,\e)$ \emph{graded quasi-isometric embedding based at}
$x\in X$, if for all $z,w \in X$:

\[ \frac{1}{\ka} d(z,w)-\e D_x(\rho;z,w) \leq d(\p (z),
\p(w) ) \leq \ka d(z,w)+\e D_x(\rho;z,w).
\]

If $F \in \mathfrak{U}$ we let $p:F\rt \Omega_F'$ be a closest
point projection and define $$\p _F:F \rt X$$ by $\p _F = \p \circ
p$.

If $x\in \Omega _F$, then using Lemma 8.3(ii), $\p _F$ is a $(\ka
, \rho, 2\kappa \varepsilon)$ graded quasi-isometric embedding
based at $x$. Also note that by Theorem 8.4, $\p_F(F)$ is
contained in a neighborhood of finitely many flats since $$\Omega
_F ' \se \big(F\cap \neu\big)_{(\varepsilon, \rho)}.$$

We fix a Weyl chamber $\ap ^+ \se \ap$ based at $\mathfrak{K} \in
X$. For any $h \in H$, let $h \ap ^+ (\infty)$ be the equivalence
class of $h \ap ^+$ in $\widehat{X}$.

For two subsets $A$ and $C$ of $X$, any point $x \in X$, and a
small number $\delta >0$, we write $A \thicksim _\delta C$ if
$$\hd\big(A\cap B_x(r) \,,\, C\cap B_x(r)\big)\leq \delta r$$ for all
sufficiently large numbers $r>0$.

At this point in \cite{Es}, a detailed argument is used to show
the analogue of the lemma below (Lemma 3.14 in \cite{Es}) for the
case when $X$ is a symmetric space.

\bigskip \noindent \textbf{Lemma 8.5} \emph{Suppose $h\ap \in
\mathfrak{U}$ for some $h \in H$. There exists a constant $\lb$
depending on $\kappa$ and $X$, and some $k\in \mathfrak{K}$
depending on $h$, such that
$$\phi_{h\ap}\big(h\app \big)\thicksim _{\lb \nre}\;k \app.$$}

Eskin's proof proceeds by first showing that if $L$ is a wall of a
flat $F\in \mathfrak{U}$, then $\phi _F$ maps $L$ into a ``graded
neighborhood" of a wall $L' \se X$. (For a definition of a graded
neighborhood see below, before the proof of Lemma 8.7.) This is
shown using the Eskin-Farb quasiflats with holes theorem and the
characterization of walls of flats in $\mathfrak{U}$ as ``coarse
intersections" of flats in $\mathfrak{U}$ (see Lemma 8.3(v)). A
key ingredient for this step is Eskin's ``no turns" lemma about
quasi-isometries of Euclidean space which respect a family of
hyperplanes. (In this case the Euclidean spaces are our flats, and
the hyperplanes are the walls of the flats.)

Since Weyl chambers are defined by the set of walls that bound
them, Eskin uses the information about the images of walls to
deduce the lemma above for symmetric spaces.

Eskin's proof of the symmetric space version of Lemma 8.5 uses the
geometry of symmetric spaces mostly to supply foundational tools
for the main argument. We will replace these tools with analogues
that hold for products of symmetric spaces and Euclidean
buildings.

The first of the foundational tools needed is Lemma 8.3---even
here Eskin's proof applied to the general case. The second tool is
Theorem 8.4 which was proved in \cite{W}. The last two tools
needed are Lemmas 8.6 and 8.7 below. They are direct analogues of
Lemmas B.1 and B.7 of \cite{Es} respectively. After proving Lemmas
8.6 and 8.7, the foundation to carry out Eskin's proof for the
general space $X$ will be in place. Then Eskin's proof applies to
establish Lemma 8.5.

\bigskip \noindent \textbf{Coarse intersections of convex polyhedra.}
Any wall, $L$, in a flat $F \se X$, divides $F$ into two
components. The closure of any such component is called a
\emph{half-space}. We define a \emph{convex polyhedron} in $X$ as
an intersection of a flat, $F$, with a (possibly empty) finite
collection of half-spaces contained in $F$. Note that flats are
convex polyhedra, as are walls.

The following lemma is an analogue of Lemma B.1 in \cite{Es}. It
allows us to replace coarse intersections of flats, walls, or
convex polyhedra with a convex polyhedron.

\bigskip \noindent \textbf{Lemma 8.6} \emph{There are constants,
$\lb _ 2$ and $\lb_3$, such that if $Q_1$ and $Q_2$ are convex
polyhedra in $X$, and if $r>\lb_2 (1+d(Q_1,Q_2))$, then there is a
convex polyhedron $P \subseteq Q_1$ such that
$$\hd \big(Q_1 \cap \nd _r(Q_2)\; ,\; P\big)\leq \lb_3 r.$$}

\pru If $Q_k\se X$ is a convex polyhedron in the flat $F_k \se X$,
 and if $F_{k,\infty}\se \xy$ and $ F_{k,\mathfrak{p}} \se \xr$
 are flats such that $F_k=F_{k,\infty} \times F_{k,\mathfrak{p}}$,
  then $$Q_k =F_k \cap \bigcap
_i \Big(H_{k,\infty,i} \times F_{k,\mathfrak{p}} \Big)\; \cap \;
\bigcap _i \Big(F_{k,\infty} \times H_{k,\mathfrak{p},i}\Big),$$
where each $H_{k,\infty ,i} \se F_{k,\infty}$ and each
$H_{k,\mathfrak{p} ,i} \se F_{k,\mathfrak{p}}$ is a half-space.

Hence, if $Q_{k,\infty} \se F_{k,\infty} $ is the convex
polyhedron given by $$Q_{k,\infty} = F_{k,\infty} \cap \bigcap _i
H_{k,\infty,i}$$ and $Q_{k,\mathfrak{p}} \se F_{k,\mathfrak{p}} $
is the convex polyhedron given by
$$Q_{k,\mathfrak{p}} = F_{k,\mathfrak{p}} \cap \bigcap _i H_{k,\mathfrak{p},i},$$
then $Q_k = Q_{k,\infty} \times Q_{k,\mathfrak{p}}$.

Note that $$ \Big[ Q_{1,\infty} \cap \nd _{r/\sqrt{2}}
(Q_{2,\infty})\Big]  \times \Big[ Q_{1,\mathfrak{p}} \cap \nd
_{r/\sqrt{2}} (Q_{2,\mathfrak{p}})\Big] $$
$$ \se $$
$$ Q_1 \cap \nd _r(Q_2) $$
$$ \se $$
$$\Big[ Q_{1,\infty} \cap \nd _{r} (Q_{2,\infty})\Big] \times \Big[
Q_{1,\mathfrak{p}} \cap \nd _{r} (Q_{2,\mathfrak{p}})\Big],$$ so
we can reduce the proof of this lemma to the separate cases of
$X=X_\infty$ and $X=X_\mathfrak{p}$. The former case is Lemma B.1
of \cite{Es}. We will prove the lemma for the latter case.

Let $Q_1$ and $Q_2$ be convex polyhedron in a Euclidean building
$X_\mathfrak{p}$. Let $F\se \xr$ be an apartment (flat) containing
$Q_1$.

Define
$$P_{d(Q_1,Q_2)}=Q_1 \cap \overline{\Nd_{d(Q_1,Q_2)}(Q_2)}$$ Since
 $Q_2$ is convex, $\overline{\Nd_{d(Q_1,Q_2)}(Q_2)}$ is convex as
  well (\cite{Br-H} Cor. II.2.5(1)). Therefore
 $P_{d(Q_1,Q_2)}$ is convex. In fact, $P_{d(Q_1,Q_2)}$
 is a convex polyhedron. Indeed, if $\mathfrak{c} \se F$ is a
 chamber, let $$
 \rho_{F,\mathfrak{c}} :\xr \rt F$$ be the retraction corresponding to $F$
 and $\mathfrak{c}$. Then $d(x,y)=d(x, \rho_{F,\mathfrak{c}}(y))$ for all
 $x \in \mathfrak{c}$ and all $y \in Q_2$. (For a good reference for
 retractions,
 and for buildings in general, see \cite{Br}.) Therefore, points in $\partial P_{d(Q_1,Q_2)}$ are determined
 by translating the region $ \rho_{F,\mathfrak{c}}(Q_2)$ a distance of
 $d(Q_1,Q_2)$. Hence, $P_{d(Q_1,Q_2)}$ is bounded by walls which
 are translates of the walls bounding $
 \rho_{F,\mathfrak{c}}(Q_2)$. Since $P_{d(Q_1,Q_2)}$ is convex,
 and since there are finitely many parallel families of walls in
 $F$, $P_{d(Q_1,Q_2)}$ is bounded by finitely many walls.

We let each $H_i \se F$ be a half-space such that
 $$P_{d(Q_1,Q_2)}=F \cap \bigcap_i H_i.$$

For any number $r\geq 0$, let $H_i^{r+} \se F$ be the half-space
that contains $H_i$, and with the additional property that
$$\hd \big(H_i \, ,\, H_i^{r+}\big) = r+d(Q_1,Q_2).$$ Define the convex polyhedron $P_r^+$ by $$P_r^+  = Q_1 \cap
\Big( \bigcap _i H_i^{r+} \Big).$$

We claim that if $r \geq 0$, then $$Q_1 \cap \nd_r (Q_2) \se
P_r^+.$$ That is, we want to prove that $$Q_1 \cap \nd_r (Q_2) \se
H_i^{r+}$$ for all $i$. To this end, let $\mathfrak{c}_i \se F$ be
a chamber that is separated from $P_{d(Q_1,Q_2)}$ by $\partial
H_i^{r+}$. Let $$
 \rho_{F,\mathfrak{c}_i} :\xr \rt F$$ be the retraction
 corresponding to $\mathfrak{c}_i$ and $F$. Since
 $\rho_{F,\mathfrak{c}_i}$is distance nonincreasing, we have that
 $$d\big(\rho_{F,\mathfrak{c}_i}
 (P_{d(Q_1,Q_2)})\,,\,\rho_{F,\mathfrak{c}_i}(Q_2)\big) \leq
 d\big(P_{d(Q_1,Q_2)}\,,\,Q_2\big)=d\big(Q_1\,,\,Q_2).$$
 Therefore, if $x\in Q_2$:
 \begin{align*}
 d\big(\partial H_i^{r+}\,,\,x\big) & \geq d\big(\partial
 H_i^{r+}\,,\, \rho_{F,\mathfrak{c}_i}(x)\big) \\
  & \geq d\big(\partial H_i^{r+}\,,\,P_{d(Q_1,Q_2)} \big) -
  d\big(P_{d(Q_1,Q_2)}\,,\,\rho_{F,\mathfrak{c}_i}(x)\big) \\
  & = r+d(Q_1,Q_2) -d\big(\rho_{F,\mathfrak{c}_i}
 (P_{d(Q_1,Q_2)})\,,\,\rho_{F,\mathfrak{c}_i}(x)\big) \\
  & \geq r.
  \end{align*}
Hence, $$Q_1 \cap \nd_r(Q_2) \se H_i^{r+}$$ as desired.

We have shown that $Q_1 \cap \nd _r (Q_2)$ is contained in a
convex polyhedron created by pushing out the walls of
$P_{d(Q_1,Q_2)}$ by a uniform distance that is linear in $r$. Next
we observe that $Q_1 \cap \nd _r (Q_2)$ also contains a convex
polyhedron created by pushing out the walls of $P_{d(Q_1,Q_2)}$ by
a uniform distance that is linear in $r$.

Indeed, since there are only finitely many walls in any flat $F'$
up to translation, there exists a positive constant $\beta <1$
depending only on $X$, such that if $Q \se F'$ is a convex
polyhedron, $s\geq 0$, and $Q(s) \se F'$ is the convex polyhedron
obtained by pushing out the walls that bound $Q$ by a distance of
$\beta s$, then
$$Q(s) \se \nd_s(Q) \cap F'.$$

Thus for any number $r\geq d(Q_1,Q_2)$, and for the set of
half-spaces $\{\,H_i\,\}$ that define $P_{d(Q_1,Q_2)}$, we let
$H_i^{r-} \se F$ be the half-space containing $H_i$ and such that
$$\hd \big(H_i \, ,\, H_i^{r-}\big) = \beta \big(r-d(Q_1,Q_2)\big).$$
And we define the convex polyhedron $P_r^-$ by $$P_r^-  = Q_1 \cap
\Big( \bigcap _i H_i^{r-} \Big),$$ so that \begin{align*} P_r^- &
\se \nd _{(r-d(Q_1,Q_2))} (P_{d(Q_1,Q_2)}) \cap F \\ & \se Q_1
\cap \nd_r(Q_2)
\end{align*}

In summary, we have shown that for $r \geq d(Q_1,Q_2)$
$$P_{r}^- \se Q_1 \cap \nd _r (Q_2) \se P_{r} ^+$$
The lemma follows since there clearly exists a constant $\lb '$
depending only on $\xr$ such that \begin{align*} \hd\big(P_{r}^-
\, ,\, P_{r}^+\big) & < \lb' [r+d(Q_1,Q_2)-\beta(r-d(Q_1,Q_2))]
\\
 & <\lb ' [r+2d(Q_1,Q_2)] \\
 & \leq \lb ' [3r].
 \end{align*}

\ep

\bigskip \noindent \textbf{Graded equivalence implies Hausdorff equivalence for Weyl chambers.} Let $\ap _\alpha \se \ap$ be a wall containing $\mathfrak{K}$.
For any collection of such walls $\{\,\ap_\alpha\,\}_{\alpha \in
\sigma}$, let
$$\ap ^+_{\sigma}=\app \cap \bigcap_{\alpha \in \sigma}\ap
_{\alpha}.$$

For any set $A \subseteq X$ and any $t>0$, we define the
\emph{graded $t$-neighborhood of} $A$ as the set
$$A[t]=\{\, x\in X \mid \mbox{there is an } a\in A \mbox{ with }
d(x,a)< td(x,\mathfrak{K}) \,\}.$$

The following lemma is a generalization of Lemma B.7 in \cite{Es}.

\bigskip \noindent \textbf{Lemma 8.7} \emph{Assume there are
three group elements $h,h_1,h_2\in H$ and that, outside of some
metric ball,
$$h \ap ^+_{\sigma} \se h_1\ap ^+[\lb \nre] \cap h_2\ap ^+[\lb
\nre].$$ If $k_1, k_2 \in \mathfrak{K}$ satisfy the condition
$$\Hd(h_i \ap^ + \,,\, k_i \ap ^+)<\infty,$$ then $$k_1 \ap
^+_{\sigma}=k_2 \ap ^+_{\sigma}.$$}

\pru A Weyl chamber $\mathfrak{C} \se X$ is a product of Weyl
chambers $\mathfrak{C}_\infty \se \xy$ and
$\mathfrak{C}_\mathfrak{p} \se \xr$. Note that
$\mathfrak{C}_\infty \times \mathfrak{C}_\mathfrak{p} \se
(\mathfrak{C}_\infty' \times \mathfrak{C}_\mathfrak{p}')[t]$
implies that, outside of a ball, $\mathfrak{C}_\infty \se
\mathfrak{C}_\infty'[t']$ and $\mathfrak{C}_\infty \se
\mathfrak{C}_\infty'[t']$ for $t'>t$. Hence, we only need to show
the case of a building since symmetric spaces are covered by Lemma
B.7 of \cite{Es}.

We can replace $h \app _{\sigma}$ by $k \app _{\sigma}$ for some
$k \in \mathfrak{K} $ such that $\Hd(h\app _{\sigma} , k \app
_{\sigma})< \infty$. Then
$$ k \app _{\sigma} \se k_1 \app [\lb \nre] \cap k_2 \app [\lb
\nre]$$ outside of a large ball.

For any $r>0$, let $a_r \in \app _{\sigma}$ be such that
$d(a_r,\app _{\alpha})>r$ for all $\alpha \notin \sigma$. By the
preceding inclusion, there exist points $a_1, a_2 \in \app$ such
that $d(ka_r,k_ia_i)\leq \lb \nre r$ for all sufficiently large
numbers $r$. Therefore, $d(k_1a_1,k_2a_2) \leq 2 \lb \sqrt[n]{\e}
r$.

There is an apartment $\ap '\se \xr$ such that, outside of a ball,
$k_i \app _\sigma \se \ap ' $ for $i=1,2$. If $k_1 \app _\sigma
\neq k_2 \app _\sigma$, then for all sufficiently large $r$, we
have $k_ia_i\in \ap ' \cap k_i \app _\sigma$ and $d(k_1a_1,k_2a_2)
> \alpha r$ for some constant $\alpha $ depending only on $\xr$.
This is a contradiction.

\ep

The proof of Lemma 8.5 only requires the case of Lemma 8.7 for
$\sigma = \emptyset$. However, the full form of Lemma 8.7 is
needed for the construction of $\pp$

\bigskip \noindent \textbf{The a.e. defined boundary function.}
Let $N<H$ be the normalizer of $A<H$. Let $\mathcal{B} (G)$
 be the Tits building for $X$. We define $\pu$ as the
simplicial subcomplex of $\mathcal{B} (G)$ given by
$$\pu=\bigcup_{h \in \mathcal{U}} \bigcup_{n\in N} hn \app(\infty).$$

We are prepared to define $$\partial \phi :\pu \rt \mathcal{B}
(G)$$ using Lemma 8.5. We let $\pp (h \app (\infty))=k
\app(\infty) $ where $k \in \mathfrak{K}$ is such that
\newline
$\phi_{h\ap}(h\app ) \thicksim _{\lb \nre}(k \app )$.

That $\pp$ is well-defined, and restricts to an isomorphism of
$\pu$ onto its image, follows from Step 4 of \cite{Es} using our
Lemma 8.7 in place of Lemma B.7 in \cite{Es}.

\bigskip \noindent \textbf{Flats are preserved.} In Section 6,
we complete $\pp$ to an automorphism of $\mathcal{B}(G)$. In Lemma
6.7, we use that apartments in $\mathcal{B}(G)$ that are contained
in $\pu$, are mapped to apartments by $\pp$. This is the content
of the lemma below. The proof is from Proposition 3.3 \cite{Es},
but we include it here as it is brief.

\bigskip \noindent \textbf{Lemma 8.8} \emph{If $F \in \mathfrak{U}$, then there is a flat
$F' \se X$ such that $\phi_{F}(F)\se \Nd_N(F')$.}

\pru Corresponding to $\p _F(F) \se X$ there is a finite set
$\mathcal{L}(\p _F)\se \widehat{X}$ of limit points (see
\cite{W}). Intuitively $\mathcal{L}(\p_F)$ is a set of equivalence
classes for finitely many Weyl chambers
$\mathfrak{C}_1,...\mathfrak{C}_k \se X$ such that
$$\hd \big(\p_F(F) \,,\, \cup_i \mathfrak{C}_i\big) < \infty .$$

Choose $x,y \in {\mathcal{L}}(\phi_{F})$ that are opposite
chambers in $\mathcal{B} (G)$. (That such chambers exist is shown
in \cite{W}.) Since $\pp$ preserves incidence relations, $\pp$ is
Tits distance nonincreasing. Therefore, $\pp ^{-1} (x)$ and $\pp
^{-1} (y)$ are opposite.

Any chamber $c \subset F (\infty )$ is contained in a minimal
gallery between   $\pp ^{-1} (x)$ and $\pp ^{-1} (y)$. Hence, $\pp
(c)$ is contained in a minimal gallery from $x$ to $y$. That is,
$\pp (c)$ is a chamber in the unique apartment containing $x$ and
$y$. Now let $F' \se X$ be the unique flat such that $F(\infty)$
contains $x$ and $y$.

\ep

\bigskip \noindent \textbf{Countable subcomplexes.}
 In Section 6 we use the following lemma to find a ``global
 sub-building'' of $\mathcal{B}(G)$ contained in $\pu$.

\bigskip \noindent \textbf{Lemma 8.9} \emph{If $V$ is a countable
 collection of chambers
in $\pu$, then there is some $h\in \tg$ such that $V\se h \pu$.}

\pru For each number $i \in \mathbb{N}$, we choose a chamber $c_i
\subset \mathcal{B} (G)$ such that $V=\{c_i\}_{i=1}^{\infty}$.
Define the set $${\mathcal{U}} _i=\{\,g\in \tg \mid gc_i \se \pu
\}.$$ Note that ${\mathcal{U}}_i \se \tg$ is co-null, so $\cap
_{i=1}^{\infty} {\mathcal{U}}_i$ is co-null. Hence, there exists
some $h^{-1} \in \cap _{i=1}^{\infty} {\mathcal{U}}_i$, and $h$
satisfies the lemma.

\ep

\bigskip \noindent \Large \textbf{9. Continuity of the boundary
function on neighborhoods of faces} \normalsize \bigskip

To complete $\pp$ to an automorphism of $\mathcal{B}(G)$ in
Section 6, we use that $\pp$ restricts to a continuous map on
simplicial neighborhoods of $(n-2)$-dimensional simplices.
Precisely, we use Lemma 9.4 below.

As with Lemma 8.5 in the previous section, our Lemma 9.4 follows
from the proof of the analogous Lemma 5.3 in \cite{Es} once a few
foundational lemmas are provided for products of symmetric spaces
and Euclidean buildings. What we require are replacements for
Lemmas B.4, B.6, and B.8 in \cite{Es}. Their analogues are listed
below as Lemmas 9.1, 9.2, and 9.3 respectively.

Recall that we defined a metric on $\widehat{X}$ in the early
portion of Section 6. We can assume that the metric is invariant
under the action of $\mathfrak{K}$. Equivalently, we assume that
the basepoint used to define the metric $\widehat{d}$ is the coset
$\mathfrak{K}\in H/\mathfrak{K} \se X$.

\bigskip \noindent \textbf{Lemma 9.1} \emph{There are constants $\nu _1$, $\nu_2$, and $\nu_3$
 depending on $X$, such that if $k_i \in \mathfrak{K}$, $z_i \in k_i \app $ with
$d(z_1,z_2)\leq \nu _1 r$, and $d(z_i , k_i \partial \app )\geq
\nu_2 r$ where $r$ is sufficiently large, then $$\widehat{d}(k_1
\app (\infty), k_2 \app (\infty))\leq \exp(-\nu_3 r).$$}

\pru The hypotheses imply the analogous hypotheses on each factor,
$\xy$ and $\xr$. On the symmetric space factor the result is
implied by Lemma B.4 of \cite{Es}, and since we have endowed
$\widehat{X}$ with the box metric, the result follows once we
establish the lemma for the case that $X$ is a Euclidean building.

Supposing $\xr$ is a Euclidean building, we let $\nu _1 =1/2$ and
$\nu _2 =1$. For the Weyl chamber $\app \se X$, we let $\alpha >1$
be the constant such that the basepoints of the sectors $\app$ and
$\app - \nd _r (\partial \app)$ are at distance $\alpha r$ from
each other for all $r>0$. Clearly $\alpha$ depends only on $\xr$.
We let $\nu _3 =\alpha /2$.

We can assume, by repositioning the direction of the geodesic rays
used to define $\widehat{d}$, that $\gamma_{\app}$ contains the
point that the sector $\app - \nd _r (\partial \app)$ is based at.
Indeed, our choice that $\gamma _{\app}(\infty) \in \app (\infty)$
is the center of mass was completely arbitrary and any point in
the interior of $\app (\infty )$ would suffice.

Now we proceed by forcing a contradiction. That is we assume that
$\widehat{d}(k_1 \app (\infty), k_2 \app (\infty))> \exp(-\alpha
r/2)$. Then $\gamma_{k_1 \app} \cap \gamma_{k_2 \app} $ is a
geodesic segment with distinct endpoints $\mathfrak{K}, x \in
\xr$, that satisfy the inequality $d(\mathfrak{K},x)<\alpha r/2$.

Let $W_x \se k_1 \ap$ be a wall containing $x$ and such that the
closure of the component of $k_1 \ap -W_x$ containing
$\mathfrak{K}$ also contains  $k_1 \app \cap k_2 \app$. Note that
the point $z_1 \in k_1 \app - \nd _r (k_1 \partial \app)$ is in
the opposite component of $k_1 \ap -W_x$ by our choice of
$\alpha$. Also by our choice of $\alpha$, $$d(z_1 , W_x) > r/2.$$

If $\mathfrak{c} \se k_1 \app$ is a chamber containing $x$, but
not contained in $k_2 \app$, then the retraction $$\rho_{k_1 \ap,
\mathfrak{c}}:\xr \rt k_1 \ap$$ corresponding to the apartment
$k_1 \ap$ and to the chamber $\mathfrak{c}$, maps $z_2$ to the
component of $k_1 \ap -W_x$ containing $\mathfrak{K}$.

Therefore, the geodesic segment from $z_1$ to $\rho _{k_1 \ap,
\mathfrak{c}}(z_2)$ passes through $W_x$. Hence, \begin{align*}
d(z_1,z_2) & \geq d(z_1, \rho _{k_1 \ap, \mathfrak{c}}(z_2)) \\
 & \geq d(z_1, W_x) \\
 & > r/2.
 \end{align*}
 This completes our contradiction.

\ep

\bigskip \noindent \textbf{Lemma 9.2} \emph{There is a constant
$\nu _4$ depending on $X$ such that for sufficiently large numbers
$Q$ and any $k_1,k_2\in \mathfrak{K}$, there are $z_i \in k_i \app
$ satisfying:} \begin{quote} $(i)$ $d(z_1,z_2)\leq Q$ \\
$(ii)$ $d(z_i,e)\leq \nu_4 \big|\log \big(\widehat{d}(k_1 \app
(\infty), k_2 \app (\infty) \big)\big|$ , and
\\
$(iii)$ $d(z_i , k_i
\partial \app )\geq \nu_5 \big | \log \big(\widehat{d}(k_1 \app (\infty) , k_2
\app (\infty) \big)\big|$ \end{quote} \emph{for some constant $\nu
_5$ which depends on $Q$ and on $X$.}

\pru Again we prove the lemma for the case $X =\xr$. The case $X
=\xy$ is Lemma B.6 of \cite{Es}, and the Lemma 9.2 follows from
the lemmas for each case.

If $\xr$ is a Euclidean building, and if $k_1 \app \cap k_2 \app$
does not contain a chamber of $\xr$, then choose $z_1\in \gamma
_{k_1 \app}$ and $z_2\in \gamma _{k_2 \app}$ to be distance $1$
away from $\mathfrak{K}$. Then the conclusion of the lemma is
satisfied for all $Q>0$ by $\nu _4$=1 and some $\nu _5$ which
depends only on the angle between $\gamma _\app$ and $\partial
\app$.

If $k_1 \app \cap k_2 \app$ does contain a chamber of $X$, then
let $z_1=z_2 \in k_1 \app \cap k_2 \app$ be the endpoint of
$\gamma _{k_1 \app } \cap \gamma _{k_2 \app }$. Now the lemma
holds for any $Q>0$, $\nu _4 =1$, and some $\nu _5$ that depends
only on the angle between $\gamma _{ \app }$ and $ \partial \app
$.

\ep

\bigskip \noindent \textbf{Lemma 9.3} \emph{Let $x,y \in X$. For any
Weyl chamber $\mathfrak{C}_x \se X$ based at $x$, there is a Weyl
chamber $\mathfrak{C}_y \se X$ based at $y$ such that
$$\hd(\mathfrak{C}_x\,,\,\mathfrak{C}_y)<\lambda ' d(x,y)$$ for
some constant $\lambda '$.}

\pru The lemma follows from Lemma B.8 of \cite{Es}, and from Lemma
4.3 of \cite{W}.

\ep

Recall that $n$ is the rank of $X$ and that for any
$(n-2)$-dimensional simplex $f \subset \pu$, we defined
$\mathcal{N}_U(f)$ as the set of all chambers in $\pu$ that
contain $f$.

We can apply the proof of Lemma 5.3 in \cite{Es} by replacing
Lemmas B.4, B.6, and B.8 of \cite{Es} with the three lemmas above
to show:

\bigskip \noindent \textbf{Lemma 9.4} \emph{If $f \subset \pu$ is a simplex of
dimension $n-2$, then $\pp |_{{\mathcal{N}}_U(f)}$ is continuous
in the Furstenberg metric.}

\bigskip

Note that Lemma 5.3 of \cite{Es} claims that $\pp
|_{{\mathcal{N}_U}(c)}$ is bi-H\"{o}lder. We only require $\pp
|_{{\mathcal{N}_U}(c)}$ to be continuous as our method for
completing $\pp$ is more algebraic, and less topological, than
Eskin's.

The condition that chambers share a wall in the above lemma is
needed so that two Weyl chambers can be simultaneously slid along
a common wall until they are based at points in $\neu$---the set
our quasi-isometry is defined on. The sliding technique does not
change their Furstenberg distance.

\bigskip

\address{
\noindent Kevin Wortman \\
Department of Mathematics \\
Cornell University \\
Malott Hall \\
Ithaca, NY 14853 \\
Email: wortman@math.cornell.edu}


\begin{thebibliography}{T-W-W-W}

\bibitem[Abl]{Abl} Abels, H., \emph{Finiteness properties of certain arithmetic groups in the function field case.} Israel J. Math., {\bf 76} (1991), 113-128.
\bibitem[Abr]{Abr} Abramenko, P., \emph{Finiteness properties of Chevalley groups over $F\sb q[t]$.} Israel J. Math., {\bf 87} (1994), 203-223.
\bibitem[Be]{Be} Behr, H., \emph{${\rm SL}\sb{3}( F\sb{q}[t])$ is not finitely presentable.}
 Proc. Sympos. ``Homological group theory''
  (Durham 1977). London Math. Soc., Lecture Notes Ser., {\bf 36}, 213-224.
\bibitem[Bl]{Bl} Blume, F., \emph{Ergodic theory.} in
\emph{Handbook of measure theory} Vol. II. North-Holland,
Amsterdam (2002), 1185-1235.
\bibitem[Bo 1]{Bo} Borel, A.,
\emph{Density and maximality of arithmetic subgroups.} J. Reine
Angew. Math., \textbf{224} (1966), 78-89.
\bibitem[Bo 2]{Bo 2} Borel, A., \emph{Linear algebraic groups.} Graduate Texts in Mathematics, No. 126,
Springer-Verlag, New York (1991).
\bibitem[Bo-Sp]{Bo-Sp} Borel, A., and Springer, T. A.,
\emph{Rationality properties of linear algebraic groups II.}
T\^{o}hoku Math. Journ., {\bf 20} (1968), 443-497.
\bibitem[Bo-T]{Bo-T} Borel, A., and Tits, J.,
 \emph{Homomorphisms ``abstraits'' de groups alg\'{e}briques simples.} Ann. Math.,
 {\bf 97} (1973), 499-571.
\bibitem[Bri-H]{Br-H} Bridson, M., and Haefliger, A., \emph{Metric spaces
 of non-positive curvature.}
Grundlehren der Mathematischen Wissenschaften, Vol. 319.
Springer-Verlag, Heidelberg (1999).
\bibitem[Bro]{Br} Brown, K., {\it Buildings.} Springer-Verlag, New York (1989).
\bibitem[Dr]{Dr} Dru\c{t}u, C., \emph{Quasi-isometric classification of non-uniform lattices in
semisimple groups of higher rank.} Geom. Funct. Anal., {\bf 10}
(2000), 327-388.
\bibitem[Es]{Es} Eskin, A., \emph{Quasi-isometric rigidity of nonuniform lattices in higher rank symmetric spaces.} J. Amer. Math. Soc., {\bf 11} (1998), 321-361.
\bibitem[E-F 1]{E-F 1} Eskin, A., and Farb, B., \emph{Quasi-flats and rigidity in higher rank symmetric spaces.} J. Amer. Math. Soc., {\bf 10} (1997), 653-692.
\bibitem[E-F 2]{E-F 2} Eskin, A., and Farb, B., \emph{Quasi-flats and rigidity in $\mathbb{H}^2 \times \mathbb{H}^2$.} Lie Groups and Ergodic Theory, Nari's Publishing House, New Delhi (1998), 75-104.
\bibitem[Fa]{Fa} Farb, B., \emph{The quasi-isometry classification of lattices in semisimple Lie groups.} Math. Res. Letta., {\bf 4} (1997), 705-717.
\bibitem[Fa-Sch]{F-S} Farb, B., and Schwartz, R., \emph{The large-scale geometry of Hilbert modular groups.} J. Diff. Geom., {\bf 44} (1996), 435-478.
\bibitem[Har]{Har} Harder, G., \emph{\"Uber die Galoiskohomologie halbeinfacher algebraischer Gruppen, III.} J. Reine Angew. Mat., {\bf 274/275} (1975), 125-138.
\bibitem[K-L]{K-L} Kleiner, B., Leeb, B., \emph{Rigidity of quasi-isometries for symmetric spaces and Euclidean buildings.} Inst. Hautes \'{E}tudes Sci. Publ. Math., {\bf 86} (1997), 115-197.
\bibitem[L-M-R]{L-M-R} Lubotzky, A., Mozes, S., and Raghunathan, M. S., \emph{The word and Riemannian metrics on lattices of semisimple groups.} Inst. Hautes \'{E}tudes Sci. Publ. Math., {\bf 91} (2000), 5-53.
\bibitem[Mar]{Mar} Margulis, G. A., \emph{Discrete subgroups of semisimple Lie groups.} Ergebnisse der Mathematik und ihrer Grenzgebeite, Springer-Verlag, Berlin-Heidelberg-New York (1991).
\bibitem[Mo]{Mo} Mostow, G. D., \emph{Strong rigidity of locally
symmetric spaces.} Princeton University Press, Princeton (1973).
\bibitem[Pl-Ra]{Pl-Ra} Platonov, V., and Rapinchuk, A.,
\emph{Algebraic groups and number theory.} Pure and Applied
Mathematics, No. 139, Academic Press, Boston, (1994).
\bibitem[Pr 1]{Pr 0} Prasad, G., \emph{Strong rigidity of \textbf{Q}-rank
 1 lattices.}, Invent. Math., \textbf{21} (1973), 255-86.
\bibitem[Pr 2]{Pr} Prasad, G., \emph{Lattices in semisimple groups
over local fields.} in \emph{Studies in algebra and number
theory.} Academic Press, New York (1979), 285-354.
\bibitem[Pr 3]{Pr 2} Prasad, G., \emph{Strong approximation for
semi-simple groups over function fields.} Ann. of Math.,
\textbf{105} (1977), 553-572.
\bibitem[Ra]{Ra} Ratner, M., \emph{On the $p$-adic and $S$-arithmetic generalizations of Raghunathan's conjectures.} Lie Groups and ergodic theory, Narosa Publishing House, New Delhi (1998), 167-202.
\bibitem[Sch 1]{S} Schwartz, R., \emph{The quasi-isometry classification of rank one lattices.} Inst. Hautes \'{E}tudes Sci. Publ. Math., {\bf 82} (1995),
133-168.
\bibitem[Sch 2]{S 2} Schwartz, R., \emph{Quasi-isometric rigidity and Diophantine approximation.} Acta Math., {\bf 177} (1996),
75-112.
\bibitem[Ta]{Ta} Taback, J., {\it Quasi-isometric rigidity for $PSL_{2}(\mathbb{Z}[1/p])$.} Duke Math. J., {\bf 101} (2000), 335-357.

\bibitem[Ti 1]{Ti 1} Tits, J., \emph{Algebraic and abstract simple groups.}  Ann. of Math.,
 \textbf{80} (1964), 313-329.
\bibitem[Ti 2]{Ti} Tits, J., \emph{Buildings of spherical type and finite BN-pairs.} Lecture Notes in Math., vol. 386, Springer-Verlag, New York (1974).
\bibitem[Ve]{Ve} Venkataramana, T. N., \emph{On superrigidity and arithmeticity of lattices in semisimple groups over local fields of arbitrary characteristic.} Invent. Math., {\bf 92} (1988), 255-306.
\bibitem[W1]{W} Wortman, K., \emph{Quasiflats with holes in reductive groups.} Preprint.

\bibitem[W2]{T-W-W} Wortman, K.,
{\it Quasi-isometries of $\mathbf{SL_n}(\mathbb{F}_q[t])$.} In
preparation.




\end{thebibliography}
\end{document}